%% file: paper.tex
\documentclass[sn-mathphys, Numbered]{sn-jnl}

\usepackage{bigfoot}

\usepackage{comment}
\usepackage{enumitem}

\usepackage{graphicx}

\usepackage{tabularx}
\usepackage{caption}
\usepackage{subfig}

\usepackage{amsmath}
\usepackage{amssymb}

\usepackage{todonotes}
\usepackage{marginnote}

\usepackage{dune}

\usepackage{pgfplots}
\pgfplotsset{compat=1.9,legend style={font=\small}}
\usepgfplotslibrary{groupplots,dateplot}
\usetikzlibrary{patterns,shapes.arrows,shapes.geometric,plotmarks}

\usepackage{multirow, makecell}

\input{macros}

\begin{document}


\title{A Jacobian-free Multigrid Preconditioner for Discontinuous Galerkin Methods applied to Numerical Weather Prediction}


\author[1]{\fnm{Philipp}   \sur{Birken}}\email{philipp.birken@math.lu.se}
\author[2]{\fnm{Andreas}   \sur{Dedner}}\email{A.S.Dedner@warwick.ac.uk}
\author*[1]{\fnm{Robert}    \sur{Kl\"ofkorn}}\email{robertk@math.lu.se}



\affil[1]{\orgdiv{Center for Mathematical Sciences}, \orgname{Lund University}, \orgaddress{Box 117}, \city{Lund}, \postcode{22100}, \country{Sweden}}
\affil[2]{\orgdiv{Mathematics Institute}, \orgname{University of Warwick}, \city{Coventry}, \postcode{CV4 7AL}, \country{UK}}

\abstract{\input{sections/abstract}}


\date{}

\maketitle

%

\input{sections/introduction}
\input{sections/dg}
\input{sections/multigrid}

\input{sections/test-cases}


\input{sections/conclusions}


\section*{Conflict of interest}
On behalf of all authors, the corresponding author states that there is no conflict of interest.

\bibliography{bibliography}


\begin{appendix}

\input{sections/appendix}



\end{appendix}

\end{document}

%% file: macros.tex
\newcommand{\phispace}{V_h^k}
\newcommand{\dgspace}{\phispace}
\newcommand{\fvspace}{V_h^{0}}

\newcommand{\grid}{{\mathcal{T}_h}}
\newcommand{\entity}{E}
\newcommand{\elem}{\entity}

\newcommand{\isec}{e}

\newcommand{\RRR}{\mathbb{R}}
\newcommand{\NNN}{{\mathbb{N}}}

\newcommand{\tol}{\text{TOL}}
\newcommand{\MG}{\text{\bf MG}}

\newcommand{\Fc}{F_{\rm c}}
\newcommand{\Fv}{F_{\rm v}}

\newcommand{\Source}{S}

\newcommand{\fluxF}{{H_{\rm c}}}
\newcommand{\fluxA}{{H_{\rm v}}}

\newcommand{\vect}[1]{\boldsymbol{#1}}

\newcommand{\vecu}{\vect{u}}

\newcommand{\basefct}{\vect{\psi}}

\newcommand{\sol}{\mathbf{U}}
\newcommand{\df}{\sol_h}
\newcommand{\U}{\sol}
\newcommand{\Ulow}{\vecu}

\newcommand{\oper}[1]{\mathcal{#1}}
\newcommand{\spcoper}{\oper{L}_h}

\newcommand{\vjump}[1]{[ \! [ {#1} ] \! ]_{\isec} }

\newcommand{\vaver}[1]{\{ \! \! \{ {#1} \} \! \! \}_{\isec} }

\newcommand{\dual}[1]{\langle \basefct, #1 \rangle}


%% file: sections/abstract.tex
Discontinuous Galerkin (DG) methods are promising high order discretizations for unsteady compressible flows. Here, we focus on Numerical Weather Prediction (NWP). These flows are characterized by a fine resolution in $z$-direction and low Mach numbers, making the system stiff. Thus, implicit time integration is required and for this a fast, highly parallel, low-memory iterative solver for the resulting algebraic systems. As a basic framework, we use inexact Jacobian-Free Newton-GMRES with a preconditioner. 
	
For low order finite volume discretizations, multigrid methods have been successfully applied to steady and unsteady fluid flows.
However, for high order DG methods, such solvers are currently lacking. 
This motivates our research to construct a Jacobian-free precondtioner for high order DG discretizations. The preconditioner is based on a multigrid method constructed for a low order finite volume discretization defined on a subgrid of the DG mesh. We design a computationally efficient and mass conservative mapping between the grids. As smoothers, explicit Runge-Kutta pseudo time iterations are used, which can be implemented in parallel in a Jacobian-free low-memory manner. 

We consider DG Methods for the Euler equations and for viscous flow equations in 2D, both with gravity, in a well balanced formulation. Numerical experiments in the software framework DUNE-FEM on atmospheric flow problems show the benefit of this approach.

%% file: sections/introduction.tex
\section{Introduction}

Discontinuous Galerkin methods has emerged as one of the most promising discretizations to replace second order finite volume methods for simulations of turbulent, time dependent, compressible fluid flow.
\begin{itemize}
	\item It is of high order, a property of critical importance for efficient Large-Eddy-Simulation (LES) of turbulent flows.
	\item It can handle the complex geometries necessary in real world applications.
	\item It has computational advantages on modern computer architectures, and lends itself well to parallelization
\end{itemize}

In this article, we make use of Discontinuous Galerkin methods (DG) on quadrilateral grids.
However, the techniques developed here can be used with other DG formulations as well,
for example, Discontinuous Spectral Element methods (DG-SEM) where
use is made of a tensor product extension from a one dimensional formulation or
other general DG methods on simplicial grids.


Many applications require implicit time stepping to avoid excessive time step
restrictions. This is the case when the mesh size varies considerably within the
domain, or when the phenomena of interest to occur over long time scales
compared to the fastest dynamics. The former occurs in particular when high
aspect ratio grids are employed, as is the case for wall bounded viscous flows,
but also for atmospheric flows where the gravitational force requires finer
resolution in the vertical direction. The latter is also a problem in numerical
weather prediction (NWP). In particular, the fastest processes are sound waves
travelling at the speed of sound, whereas the phenomena of interest travel at
the flow velocity. For low Mach flows, as is often the case in NWP, this is a
serious issue.

In this article, we will consider fully implicit time integration, which leads
to one large global system of equations. Complementary approaches that we will
not follow in this text, are implicit-explicit (IMEX) methods, and horizontal
explicit, vertical implicit (HEVI) ones that have seen a lot of application in
NWP \cite{hevi:15,giraldohevi:24}. These make use of a geometric or term wise
splitting of the equations and use different time integration methods for the
different parts. The technique suggested here can be employed subsequently in
those methods as well.

Implicit time stepping requires solving large sparse linear systems of
equations, which in turn requires efficient solvers to be competitive
computationally.  As pointed out in
\cite{bastianMatrixfreeMultigridBlockpreconditioners2019}, to achieve the peak
floating point performance on modern processors, it is necessary to reuse data
loaded into the CPU, i.e. the arithmetic intensity of the algorithms must be
sufficiently high.  Algorithms based on frequent matrix-vector products with
assembled matrices are unlikely to achieve sufficient arithmetic intensity, and
it is therefore necessary to consider matrix-free or Jacobian-free algorithms
\cite{dunecodegen, kronbichler2019matrixfree} and preconditioners suitable for
such solvers.  The lack of efficient and grid independent solvers for DG
discretizations of compressible time dependent flow problems limits the
applicability of high order discretizations in many practical applications.
Previous attempts are for example \cite{baghrt:09, bigahm:13, Carr2015,
frpepa:22, watrwija:22, Puppo2024}.  Robust and efficient preconditioning
techniques for Jacobian-free Newton-Krylov (JFNK) solvers are therefore an
active area of research.

Multigrid methods have been shown to scale optimally in many settings, in
particular for elliptic problems, but also for finite volume discretizations of
compressible flows \cite{caujam:01,bibuja:19,
birkenNumericalMethodsUnsteady2021}. For the multigrid method to be
Jacobian-free, the smoother must be Jacobian-free, restricting the possible
choices. Recently, a new Jacobian-free multigrid preconditioner for use in JFNK
methods was suggested in \cite{birken:19} and then extended in
\cite{kasiECCOMAS:21,bidekaklECCOMAS:24}. The approach uses a geometric
multigrid method defined on a sub cell finite volume discretization (also known
as low-order-refined (LOR) preconditioning). It is therefore similar to
$p$-multigrid, but allows to make use of existing smoothers for FV methods, such
as \cite{ birken:12, bigahm:13}.

This motivates us to build on that foundation. We follow the framework from
\cite{bbbkmz:16}, which makes use of high order time integration, and a JFNK
solver with a smart choice of tolerances. We extend the preconditioner named
above to the discretized 2D viscous compressible flow equations with a source
term for gravity, as is common for atmospheric flows. Thereby, we make use of
the technique from \cite{cosmodg:14} to obtain a well balanced scheme.

The paper is organized as follows.  In Sections \ref{seq:discretization_spatial}
and \ref{TimeDisc} we briefly recall the main building blocks of the DG
discretization and implicit time stepping.  In Sections \ref{sec:solver_start}
and \ref{sec:multigrid} we introduce the proposed solver and preconditioner for
DG discretizations.  We then investigate in Section \ref{sec:testcases} the
effectiveness of this approach for various atmospheric test cases.


\section{Governing equations, Discretization and Solvers}

\label{sec:equations}
We consider a general class of time dependent nonlinear advection-diffusion-reaction problems
for a vector valued function $\sol\colon(0,T)\times\Omega\to\RRR^r$
with $r\in\NNN^+$ components of the form
\begin{eqnarray}
\label{eqn:general}
\label{eqn:ns}
\label{eqn:general_op}
  \partial_t \sol  = \oper{L}(\sol) &:=& - \nabla \cdot \big(\Fc(\sol)
                                         +\Fv(\sol,\nabla\sol) \big)
                                         \ \ \mbox{ in } (0,T] \times \Omega
\end{eqnarray}
in $\Omega \subset \RRR^d$, $d=1,2,3$. Suitable initial and boundary conditions have
to be added. $\Fc$ describes the convective flux,
$\Fv$ the viscous flux, and
$\Source$ a source term.
Note that all the coefficients in the partial differential equation are
allowed to depend explicitly on the spatial variable $x$ and on time $t$
but to simplify the presentation we suppress this dependency in our
notation.

For the discretization, we use a method of lines approach based on first
discretizing the differential operator in space using a DG approximation, and then solving the resulting system of ODEs using a time stepping scheme.

%% file: sections/dg.tex

\subsection{Discontinuous Galerkin Discretization}
\label{sec:dune}
\label{seq:discretization_spatial}
Given a tessellation $\grid$ of the computational domain $\Omega$ with
$\cup_{\elem \in \grid} \elem = \Omega$
we denote with $\Gamma_i$ the set of all intersections between two
elements of the grid $\grid$, and the set of all
intersections, also with the boundary of the domain $\Omega$, is denoted by $\Gamma$.


We consider a discrete space $\phispace$ spanned by
Lagrange type basis functions $\basefct_i(x)$ based on tensor product Gauss-Legendre (GL) quadrature nodes.
This yields diagonal mass matrices as discussed in detail
in \cite{dunefemdg:21, kopriva:10, kopriva:02}.

After fixing the grid and the discrete space,
we seek
\[\df(t,x) = \sum_i \sol_i(t)\basefct_i(x) \in \phispace,\]
by discretizing the spatial operator
$\oper{L}(\sol)$ in \eqref{eqn:general_op}
with boundary conditions
by defining for all test functions $\basefct \in \phispace$,
\begin{equation}
\label{convDiscr}
\dual { \spcoper(\df) } := \dual{ K_h(\df) } + \dual{ {I}_h(\df) }.
\end{equation}
Hereby, we have the element integrals
\begin{eqnarray}
\label{eqn:elementint}
   \dual{ {K}_h(\df) } &:=&
      \sum_{\elem \in \grid} \int_{\elem}
      \big( (\Fc(\df) - \Fv(\df, \nabla \df ) ): \nabla\basefct + \Source(\df)
      \cdot \basefct \big),
\end{eqnarray}
and the surface integrals (by introducing appropriate numerical fluxes
$\fluxF$, $\fluxA$ for the convection and diffusion terms, respectively)
\begin{eqnarray}
\label{eqn:surfaceint}
   \dual{ {I}_h(\df) } &:=&
      \sum_{e \in \Gamma_i} \int_e \big(
      \vaver{\Fv(\df, \vjump{\df} )^T : \nabla\basefct} +
      \vaver{\Fv(\df, \nabla\df)} : \vjump{\basefct} \big) \nonumber \\ &&
    - \sum_{e \in \Gamma} \int_e  \big(\fluxF(\df) - \fluxA(\df,\nabla\df)\big):
      \vjump{\basefct},
\end{eqnarray}
with
 $\vaver{ \sol }$,
$\vjump{ \sol }$ denoting the classic
average and
jump of $\sol$ over $e\in\Gamma$,
respectively.
The convective numerical flux $\fluxF$ is chosen to be the HLLC flux (see \cite{Toro99} for details).

A wide range of diffusion fluxes $\fluxA$ can be found in the literature and many of these fluxes are available in \dunefemdg, for example, Interior Penalty and variants, Local DG, Compact DG 1 and 2 as well as Bassi-Rebay 1 and 2 (cf. \cite{cdg2:12,dunefemdg:17}).

We then require the Galerkin condition
\begin{eqnarray}
\label{eqn:dgform}
 \frac{d}{dt}\int_\Omega\df\cdot\basefct \,dx = \dual { \spcoper(\df) }
\end{eqnarray}
for all test functions $\basefct \in \phispace$.
\subsection{Temporal discretization}
\label{TimeDisc}

%
Equation \eqref{eqn:dgform} is a system of ODEs for the coefficients of $\sol(t)$ which reads
\begin{eqnarray}
  \label{eqn:ode}
  \sol'(t) &=& f(\sol(t))  \mbox{ in } (0,T]
\end{eqnarray}
with $f(\sol(t)) = M^{-1}\spcoper(\df(t))$ where the mass matrix $M$ is given element wise by
\[m_{ij}=\int_{E} \varphi_i\varphi_j =
  \omega_i\delta_{ij}.\]
The initial value $\sol(0)$ is given by the $L_2$ projection of $\sol_0$ onto $\phispace$.
%

We use Singly-Diagonally Implicit Runge Kutta (SDIRK) methods because of
their simplicity and relatively low cost compared to fully implicit schemes. A SDIRK method is defined by a lower triangular Butcher tableau with constant nonzero diagonal elements. In every stage, a nonlinear equation system of the form
\begin{equation}
\label{eqn:newton}
\mathbf{G}(\mathbf{U}_s) := \mathbf{U}_s - \alpha \Delta t f(\mathbf{U}_s) - \bar{\mathbf{U}}_s = \mathbf{0}
\end{equation}
is to be solved, where $U_s$ is the stage value at stage $s$. $\bar{U}_s$ and $\alpha$ are the previous stage value and the diagonal entry in the Butcher tableau.

\begin{wrapfigure}{r}{0.25\textwidth}
\begin{equation}
\renewcommand\arraystretch{1.2}
\begin{array}
{c|cc}
\alpha & \alpha \\
1      & 1-\alpha & \alpha \\
\hline
       & 1 - \alpha & \alpha
\end{array}
\label{eqn:ellsiep}
\end{equation}
\end{wrapfigure}
Specifically, we use a method of order 2 with 2 stages known as Ellsiepen's method \eqref{eqn:ellsiep}, or SDIRK2,
with $\alpha = 1 - \frac{\sqrt{2}}{2}$, see \cite{Ellsiepen:99}.

%% file: sections/multigrid.tex
\input{sections/multigrid-solver.tex}

\input{sections/multigrid-loap.tex}

\input{sections/multigrid-transfer.tex}

\input{sections/massconserve.tex}

\input{sections/multigrid-fvapprox.tex}

%% file: sections/multigrid-solver.tex
\subsection{Solver}
\label{sec:solver_start}
The systems~\eqref{eqn:newton} are solved using a preconditioned Jacobian-free Newton-Krylov
method (see \cite{knoll:04}). The preconditioning strategy is the main topic of this article. We use a multigrid preconditioner based on a low order discretization of the problem.

Newton's method is applied to solve \eqref{eqn:newton} with initial guess $\mathbf{U}_s^{0} = \bar{\mathbf{U}}_s$:
\begin{eqnarray}
\label{eqn:krylov}
  \mathbf{G}'(\mathbf{U}_s^k) \delta \mathbf{U}_s &=& -\mathbf{G}(\mathbf{U}_s^k), \\
  \mathbf{U}_s^{k+1} &=& \mathbf{U}_s^k + \delta \mathbf{U}_s. 
\end{eqnarray}
The iteration is terminated when
\begin{equation}
    |\mathbf{G}(\mathbf{U}_s^k)| < \tol \ |\mathbf{G}(\mathbf{U}_s^0)|
\end{equation}
where $\tol$ is a specified tolerance. The linear system \eqref{eqn:krylov} is solved 
using GMRES with a preconditioner based on a multigrid method, which will be
explained in the following subsections.  
The tolerance for the termination of the Krylov method is selected 
adaptively using the second Eisenstat-Walker
criterion~\cite{EisenstatWalker:96} (setting the parameters $\gamma = 0.1$ and $\alpha = 1$).

GMRES is used in a Jacobian-free manner, i.e. using a finite
difference approximation of the Jacobian-Vector product, to avoid the storage of
large Jacobian matrices, which for higher order DG methods are prohibitively large:
\begin{equation}
    \mathbf{G}'(\mathbf{U}) \mathbf{y} \approx \frac{1}{\epsilon} \bigg[\mathbf{G}(\mathbf{U} + \epsilon {\mathbf y}) - \mathbf{G}(\mathbf{U})\bigg],
\end{equation}
where $\epsilon = \sqrt{\epsilon_{mach}} / |\mathbf{y}|$ and $\epsilon_{mach}$ is the machine precision.

%% file: sections/multigrid-loap.tex
\subsection{Low order approximation preconditioner}
\label{section:preconditioner}

The main difficulty here now lies in finding an effective Jacobian-free preconditioner.
The setup excludes some common options such as Gauss-Seidel or incomplete LU preconditioners.
Instead we choose a geometric agglomeration multigrid preconditioner based on an auxiliary low order discretization of the problem, described in \cite{birken:19,kasiECCOMAS:21}.

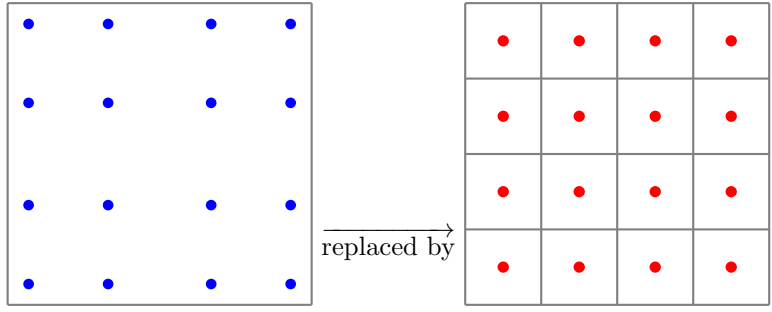
\begin{figure}
        \begin{center}
  \begin{tikzpicture}[every node/.style={inner sep=0,outer sep=0},scale=4]
  \newcommand*{\xMin}{0}%
  \newcommand*{\xMax}{1}%
  \newcommand*{\yMin}{0}%
  \newcommand*{\yMax}{1}%
  %
  \newcommand*{\xMinPt}{0.069}%
  \newcommand*{\xMaxPt}{0.862}%
  \newcommand*{\yMinPt}{0.069}%
  \newcommand*{\yMaxPt}{0.931}%
  \newcommand*{\xMidPt}{0.339}%

  \foreach \i in {0,4} {
    \draw [thick,gray] (\xMax/4 * \i,\yMin) -- (\xMax/4 * \i,\yMax);
  }
  \foreach \i in {0,4} {
    \draw [thick,gray] (\xMin,\yMax/4 * \i) -- (\xMax,\yMax/4 * \i);
  }
  \foreach \y in {-1+\xMinPt, -\xMidPt, \xMidPt, \xMaxPt}
  \foreach \x in {-1+\xMinPt, -\xMidPt, \xMidPt, \xMaxPt} {
    \draw node [circle,fill=blue,scale=4] at (\xMax/2*\x+\xMax/2, \yMax/2*\y+\yMax/2) {};
  }
  \end{tikzpicture}
  $\begin{aligned}\overrightarrow{\text{replaced by}} \\ \\[25pt]\end{aligned}$
  \begin{tikzpicture}[every node/.style={inner sep=0,outer sep=0},scale=4]
  \newcommand*{\xMin}{0}%
  \newcommand*{\xMax}{1}%
  \newcommand*{\yMin}{0}%
  \newcommand*{\yMax}{1}%
  \foreach \i in {0,4} {
    \draw [thick, gray] (\xMax/4 * \i,\yMin) -- (\xMax/4 * \i,\yMax);
  }
  \foreach \i in {0,4} {
    \draw [thick,gray] (\xMin,\yMax/4 * \i) -- (\xMax,\yMax/4 * \i);
  }

  \foreach \i in {1,...,3} {
    \draw [thick,gray] (\xMax/4 * \i,\yMin) -- (\xMax/4 * \i,\yMax);
  }
  \foreach \i in {1,...,3} {
    \draw [thick,gray] (\xMin,\yMax/4 * \i) -- (\xMax,\yMax/4 * \i);
  }
  \foreach \y in {-0.75, -0.25, 0.25, 0.75}
  \foreach \x in {-0.75, -0.25, 0.25, 0.75} {
    \node[red,scale=4] at (\xMax/2*\x+\xMax/2, \yMax/2*\y+\yMax/2) {.};
  }
  \end{tikzpicture}
  \vspace{-0.5cm}
          \caption{(Left) Gau\ss-Legendre (GL) nodes for $k=3$
          in one grid cell and (right) FV cell averages in subgrid.
          The specific position of interpolation nodes for the DG method
          is only of interest for DG-SEM variants.}
          \label{fig:lgl_fv_points}
      \end{center}
\end{figure}

We make use of a finite volume discretization, defined by the
description in Section \ref{seq:discretization_spatial} with
polynomial degree $k=0$. It is defined on a subgrid of the DG grid in such a way that the two discretizations have the same number of degrees of freedom.
The low order spatial discretization on the refined grid defines a function
space $\fvspace$ consisting of piece-wise constant functions and an operator $f_{low}$, such that the low order spatial discretization of the problem reads
\begin{equation}
  \Ulow^{\prime}(t) = f_{low}(\Ulow(t)),
\end{equation}
where $\Ulow$ is the corresponding vector of the degrees of freedom. An
example of a subgrid using an equidistant subdivision is shown in
Figure~\ref{fig:lgl_fv_points} for a DG space of order three.

Using the same temporal discretization as for the high order problem, the spatially low order discretization gives rise to a nonlinear system corresponding to eqn. \eqref{eqn:newton}
\begin{eqnarray}
    \mathbf{g}(\Ulow_s) := \Ulow_s - \alpha \Delta t f_{low}(\Ulow_s) - \bar{\Ulow}_s =\mathbf{0},
\end{eqnarray}
and the linear system
\begin{equation}
\label{eqn:lowkrylov}
\mathbf{g}'(\Ulow_s^k) \delta \Ulow_s = -\mathbf{g}(\Ulow_s^k)
\end{equation}
corresponding to eqn. \eqref{eqn:krylov}.

Assuming we have a preconditioner $\mathbf{q}^{-1}(\Ulow_{s}^k) \approx  \mathbf{g}'(\Ulow_{s}^k)^{-1}$ for \eqref{eqn:lowkrylov}, we define a preconditioner for eqn. \eqref{eqn:krylov} to be
\begin{equation}
     \label{eqn:preconditioner}
     \mathbf{Q}^{-1}:= \mathbf{T}^{-1}\mathbf{q}^{-1}\mathbf{T},
\end{equation}
where $\mathbf{T}$ is a \textit{transfer function} between the two discretizations (see Section \ref{sec:transfer-function}).
Note that the transfer function can also be combined with a pre or post
smoothing step, which is discussed in Section~\ref{sec:multigrid}.

%% file: sections/multigrid-transfer.tex
\subsubsection{Transfer functions}
\label{sec:transfer-function}

We now define mappings between the degrees of freedom of the high and the low order discretization. 
The vectors $\U_s, \U_s^k, \bar{\U}_s$ in the original discretization represent functions in $\dgspace$.
Correspondingly, $\Ulow_{s}, \Ulow_{s}^k, \bar{\Ulow}_{s}$ 
represent solutions of the low order discretization, i.e. functions in $\fvspace$.
The transfer functions between $\dgspace$ and $\fvspace$ are defined 
such that the functions are close \textit{in some sense}.

An obvious choice for $\mathbf{T}$ would be the $L^2$ projection. However, in our experience, the simpler mapping of interpolating the polynomials in $\dgspace$ 
in the cell centers of the cells defining $\fvspace$ gives the same preconditioner performance, while being much easier to implement on unstructured grids
and requiring less computational cost. 
A detailed comparison of different approaches shows the interpolation based
transfer functions offer the best performance in terms of preconditioner
iterations and computing time \cite{kasimaster:21, kasiECCOMAS:21} and was therefore chosen for the experiments in
this paper. Further discussions on this topic can be found in
\cite{witherdenNodalPointSets2021}.

%% file: sections/massconserve.tex
\subsubsection{Conservation of mass}
For the operator in \eqref{eqn:preconditioner} to conserve mass, both the low order preconditioner and the transfer function have to do so. See \cite{birlin:22,linbir:23} for a discussion on mass conservation of iterative methods and preconditioners. If the transfer function is based on $L_2$ projection, it conserves mass, but the interpolation introduces small conservation errors.
If the problem at hand is sensitive to mass change, it is possible to add a {\it mass fix} to make the transfer function mass conservative. 

The mass fix adds a term to the interpolation transfer function that is the average of the mass difference between the high order and the low order approximation in every DG cell:
\begin{equation}
\label{eqn:transfer-mass-fix}
\mathbf{T}^{mf}_\elem\mathbf{u} := \mathbf{T}_\elem\mathbf{u} - \frac{(m_\elem(\mathbf{T}_\elem\mathbf{u}) - m_\elem(\mathbf{u}))}{|\elem|} \mathbf{1}.
\end{equation}
Here, $\mathbf{T}_\elem$ is the block of $\mathbf{T}$ relevant to the grid cell $\elem$ and $m_{\elem}(\cdot) = \int_{\elem} \cdot$ is the mass in cell $\elem$. For the low order approximation given by $\mathbf{T}_\elem\mathbf{u}$, this is straighforward to compute.

For the high order approximation $\mathbf{u}$, it is possible to avoid evaluating this function in the DG nodes as well. To this end, the evaluations at the
cell centers of the finite volume cells can be reused by using suitable quadrature weights
such that these form a quadrature of the correct order. These weights are listed
in the Appendix for $k=3$.

For the test cases in this article we observed no difference between the interpolation transfer function 
and the mass fix interpolation transfer function in terms of preconditioner efficiency or quality of the solution.

%% file: sections/multigrid-fvapprox.tex
\subsubsection{Multigrid preconditioner for finite volume approximations}
\label{sec:multigrid}

We start with a sequence of increasingly refined quad meshes $\mathcal{T}_{l}$ for $l=0\ldots L$,
where $\mathcal{T}_{l}$ is generated by splitting each quad
$\elem\in\mathcal{T}_{l-1}$ into $2^d$ new child quads denoted by $C(\elem)$.
The approximation on level $l$ is represented by a vector $\mathbf{U}_l$.

The method is defined by choosing restriction and prolongation operators, 
$\mathbf{R}, \ \mathbf{P}$, and a \textit{smoother}, $\mathbf{S}(\cdot,\cdot)$.

For the restriction $\mathbf{R}$ we use the agglomeration restriction, mapping the cell averages on the quads $C(E)$ of the fine tessellation $\mathcal{T}_l$ to the volume averaged cell average on quad $E$ of the coarser tessellation $\mathcal{T}_{l-1}$:
\begin{equation}
    \big(\mathbf{R}\mathbf{U}_l\big)_\elem = \frac{1}{|\elem|}\sum_{q \in C(\elem)} |q|U_{l_q}, \quad \forall \elem\in \mathcal{T}_{l-1}.
\end{equation}

For prolongation we use the injection operator $\mathbf{P}$. It assigns the cell average of the agglomerate to all split quads. It can be expressed as
\begin{equation}
\big(\mathbf{P} \mathbf{U}_{l}\big)_q  =  U_{l_\elem}, \quad \forall q \in C(\elem), \quad \elem\in \mathcal{T}_{l}.
\end{equation}

The smoother $\mathbf{S}(\cdot, \cdot)$ is a linear iterative method that converges to the solution of \eqref{eqn:lowkrylov}.
It is designed to quickly reduce high frequency components of the error.
We use a smoother based on pseudo time iteration, described in section \ref{sec:pseudo-time-iteration}.

The multigrid algorithm consists of a V- or W-cycle with pre- and post-smoothing.
On the coarsest grid level we apply the smoother twice.
The multigrid algorithm for solving equation \eqref{eqn:lowkrylov} can be written as
\begin{enumerate}
  \item Compute $\MG_l(\mathbf{x}, \mathbf{b}; \ \Ulow_s^k)$:
  \item $\mathbf{x} := \mathbf{S}(\mathbf{x}, \mathbf{b})$
  \item If $l > 0$ 
    \begin{itemize}
      \item $\mathbf{r} := \mathbf{R}(\mathbf{g}'(\Ulow_s^k) \mathbf{x} - \mathbf{b})$
      \item $\mathbf{v} := \mathbf{0}$
      \item Repeat once for V-cycle and twice for W-cycle 
        \begin{itemize}
          \item[] $\mathbf{v} := \MG_{l-1}(\mathbf{v}, \mathbf{r}; \ \mathbf{R}\Ulow_s^k) $
        \end{itemize}
      \item $\mathbf{x} := \mathbf{x} - \mathbf{P}\mathbf{v}$
    \end{itemize}
    \item $\mathbf{x} := \mathbf{S}(\mathbf{x}, \mathbf{b})$
    \item Return $\mathbf{x}$
    \end{enumerate}

The preconditioner is obtained by calling this routine with inital guess $\mathbf{x}=\mathbf{0}$. To make this efficient, we add a flag that avoids MatVecs with the zero vector. 

\subsubsection{Smoother / Pseudo time iteration}
\label{sec:pseudo-time-iteration}
A common approach to solve linear systems such as \eqref{eqn:lowkrylov} arising 
from computational fluid dynamics models is so called pseudo-time stepping.
Smoothing can be both done during the multigrid cycle, i.e., on the
finite-volume data before or after prolongation/restriction. Smoothing is
also possible on the DG data before the transfer to the finite-volume grid
or after reconstruction of the DG data. In both cases we use a similar pseudo-time
stepping approach which we describe here only for the finite-volume data.

Instead of solving \eqref{eqn:lowkrylov} directly, a pseudo time variable is introduced
\newcommand{\pseudo}{\mathbf{w}}
\begin{equation}
     \frac{\delta \pseudo}{\delta \tau} = -\mathbf{g}(\Ulow_s^k) - \mathbf{g}'(\Ulow_s^k) \pseudo
\end{equation}
to construct an ODE with a steady state at the solution of the linear system.
The ODE is stable if the original spatial discretization of the system is stable.
In the next step, an explicit Runge-Kutta (ERK) time integration scheme is used to 
integrate in pseudo time from an initial guess, namely the current approximation $\mathbf{x}$ in the pseudocode. 

These methods have several parameters, namely Runge-Kutta coefficients and a time step in pseudotime, $\Delta t^*$. The latter is computed based on a pseudo CFL number. This means that on coarser grids, larger time steps will be taken. These are chosen based on \cite{birken:12}. There, the parameters are optimized for certain physical CFL numbers. In the experiments, the smoother is the one stage explicit Euler scheme, which we found to be the most efficient.

%% file: sections/test-cases.tex
\section{Numerical Experiments}
\label{sec:testcases}
\input{sections/test-cases/euler-model}
\subsection{Non-hydrostatic inertia gravity}
\input{sections/test-cases/inertia-gravity}

\subsection{Rising warm air bubble}
\input{sections/test-cases/rising-bubble}
\subsection{Density current}
\input{sections/test-cases/density-current}

%% file: sections/test-cases/euler-model.tex
The governing equations for the test cases are two dimensional 
viscous compressible flow equations. Expressed as in Eq. \eqref{eqn:general} they are given by
\begin{align}
    \label{eqn:compressible-euler}
    F_c(\mathbf{U}) = \begin{bmatrix}
    \rho u & \rho w \\
    \rho u^2 + p & \rho uw \\
    \rho u w & \rho w^2 + p \\
    u\rho \theta & w\rho \theta
    \end{bmatrix}, \quad  
    F_v(\mathbf{U}, \nabla\mathbf{U}) = \mu \rho \begin{bmatrix}
    \mathbf{0} \\
    \nabla u \\
    \nabla w \\
    \nabla \theta
    \end{bmatrix}, \quad
    S_i(\mathbf{U}) = \begin{bmatrix}
    0 \\ 0 \\ -\rho g \\ 0
    \end{bmatrix},
\end{align}
with $S_e(\mathbf{U}) = \mathbf{0}$ in the test cases if nothing else is specified.
If the diffusion term is neglected, one recovers the Euler equations with a gravity source term. These will be use for several of the test cases. 

It is common for atmospheric test cases to consider a perturbation of a 
\emph{environmental atmosphere}, i.e., a stationary solution 
$\bar{\mathbf{U}} = (\bar{\rho}, \bar{\rho}\bar{\mathbf{v}}, \bar{\rho}\bar{\theta})$, which fulfill \eqref{eqn:compressible-euler} with $\mu=0$. To avoid numerical errors in balancing the environmental atmosphere from polluting the simulation results, we use the well-balancing scheme introduced in \cite{dissDedner}, which is based on applying our DG method to the system for the perturbations $\mathbf{U}'$:
\begin{equation}
\label{eqn:perturbed-euler}
\partial_t \mathbf{U}' + \nabla \cdot (F_c(\mathbf{U}'+\bar{\mathbf{U}}) - F_c(\bar{\mathbf{U}}) - F_v(\mathbf{U}'+\bar{\mathbf{U}}, \nabla(\mathbf{U}'+\bar{\mathbf{U}})) = S_i(\mathbf{U}'+\bar{\mathbf{U}}) - S_i(\bar{\mathbf{U}}).
\end{equation}

To apply DG to \eqref{eqn:perturbed-euler}, we rewrite the equation as
\begin{equation*}
    \partial_t \mathbf{U}' + \nabla \cdot (F_c^{pert}(\mathbf{U}') - F_v^{pert}(\mathbf{U}')) = S^{pert}(\mathbf{U}')
\end{equation*}
with
\begin{align*}
    &F_c^{pert}(\mathbf{U}')  = F_c(\mathbf{U}'+\bar{\mathbf{U}}) - F_c(\bar{\mathbf{U}}),
    \quad F_v^{pert}(\mathbf{U}') = F_v(\mathbf{U}'+\bar{\mathbf{U}}, \nabla(\mathbf{U}'+\bar{\mathbf{U}})), \\
    &S^{pert}_i(\mathbf{U}') = S_i(\mathbf{U}'+\bar{\mathbf{U}}) - S_i(\bar{\mathbf{U}}),
\end{align*}
and apply the DG procedure as described in Section \ref{seq:discretization_spatial}. 
The corresponding numerical fluxes $\hat{F}_c^{pert}(\mathbf{U}')$ and $\hat{F}_v^{pert}(\mathbf{U}')$ now read
\begin{equation*}
\hat{F}_c^{pert}(\mathbf{U}') = \hat{F}(\mathbf{U}'+\bar{\mathbf{U}}) - \hat{F}(\bar{\mathbf{U}}), \quad \hat{F}_v^{pert}(\mathbf{U}') = \hat{F}_v(\mathbf{U}'+\bar{\mathbf{U}}, \nabla(\mathbf{U}'+\bar{\mathbf{U}})).
\end{equation*}
Applying the DG method to \eqref{eqn:perturbed-euler} guarantees 
that $\mathbf{U}' \equiv 0$ is maintained on the discrete level, 
if $\mathbf{U}' $ vanishes initially. Furthermore, this version of 
the system has significant advantages in adaptive 
computations as shown in \cite{cosmodg:14}. 
Notice that this approach requires that the environmental atmosphere 
is given a-priori.

In all tests we compute a DG solution based on tensor product Legendre polynomials of degree $k=3$ on a given refinement level $l$ of a cartesian grid.
The multigrid preconditioner then uses $l+2$ levels since
going from DG with order $k=3$ to the finest FV grid adds two level as shown in Figure \ref{fig:lgl_fv_points}.
To test the effectiveness of the multigrid preconditioner we use an implicit SDIRK2 with a fixed time step $\Delta t$ and fixing a relatively coarse tolerance of
$\tol = 10^{-3}$ for the Newton method if not mentioned otherwise. For the linear solver we use a matrix free GMRES method where the tolerance is selected adaptively using the second Eisenstat-Walker criterion as mentioned previously.
We then compare the number of iterations of the GMRES solver for the non preconditioned and different multigrid setups.

For the preconditioner there are a variety of configuration parameters which are
encrypted in the following key for our experiments:
\begin{equation}
  \text{\bf mg} \, a \,b \, c \, d\, e\, f \, G
\end{equation}
The \textbf{mg} is an abbreviation for multigrid and the other letters have the following meaning:
{
\begin{description}[topsep=5pt, noitemsep, align=left, labelwidth=1cm,labelindent=1cm]
  \item[a,b:] Number of \textbf{pre} and \textbf{post} smoothing steps on the DG solution
  \item[c,d:] Number of \textbf{pre} and \textbf{post} smoothing steps on the finest FV level
  \item[e,f:] Number of \textbf{pre} and \textbf{post} smoothing steps on the intermediate FV levels
  \item[G:] Either \textbf{V} or \textbf{W}, denoting the cycle of the multigrid method
\end{description}
}
For example, the configuration $\text{\bf mg}111111V$ means multigrid in V-cycle mode with one pre and post smoothing step on the DG solution as well as on all FV levels.
 
To obtain a reference solution, we in addition include results for
an explicit 4-stage 3rd order SSP method described in \cite{Ketcheson:08}
and implemented in \dunefemdg~\cite{dunefemdg:21} with a time step $\Delta t$ close to the stability limit of the explicit method.

%% file: sections/test-cases/inertia-gravity.tex
This test case, introduced in \cite{skamarock-klemp:94},
considers the evolution of a potential temperature perturbation.
The perturbation is so small that the bubble does not have enough buoyancy to rise,
but rather oscillates in the vertical direction, creating a specific wave like
pattern in the horizontal direction. The domain is a tube of dimension $[0,\ 300\,000] \times [0,\ 10\,000]$km$^2$.
The given environmental atmosphere has constant Brunt-Väisälä frequency $N=10^{-2}$ s$^{-1}$,
the surface pressure $p_0=10^5$Pa, the surface temperature $T_0=250$K, and a
constant mean flow $\tilde{u}=20$m/s and $\tilde{w}=0$.

The governing equations for this test case are the Euler equations with gravity source term. The potential temperature, temperature and pressure values of the environmental atmosphere are given by
\begin{align}
    &\tilde{\theta}(x, z) = T_0 \exp\frac{z}{H}, \\
    &\tilde{T}(x, z) = T_0\big(\alpha - (\alpha - 1)\exp\frac{z}{H}\big), \\
    &\tilde{p}(x, z) = p_0 \exp\bigg( \frac{c_p}{R_d}\bigg(\ln\frac{\tilde{T}}{T_0}-\frac{z}{H}\bigg)\bigg)
\end{align}
with constants $H$ and $\alpha$ defined as
\begin{align}
    H = \frac{g}{N^2}, \quad \alpha = \frac{gH}{c_p T_0}.
\end{align}
For this test case, we choose $c_p = 1.005$J/(kg K), $c_v = 717.95$J/(kg K), and $g = 9.80665$m/s$^2$.

The initial perturbation of the potential temperature is given by
\begin{align}
    \theta'(x, z) = \theta_c \bigg(1 + \bigg(\frac{x-x_c}{a}\bigg)^2\bigg)^{-1} \sin \frac{\pi z}{Z},
\end{align}
with $Z=10^4$m preserving the pressure values. The maximal deviation from the stratified atmosphere $\theta_c = 0.01$K, the center of the perturbation bubble on the x-axis $x_c=100,000$m, and $a=5,000$m. Introducing the perturbation $\theta'$ of the potential temperature into the environmental atmosphere while preserving pressure finally defines the initial state $\mathbf{U} = (\rho, \rho\mathbf{v}, \rho \theta)$ through
\begin{align}
    \label{eq:add-perturbation}
    \rho = \frac{p_0^{R_d/c_p} \tilde{p}^\frac{1}{\gamma}}{R_d (\tilde{\theta} + \theta')}, \quad \mathbf{v} = \tilde{\mathbf{v}}, \quad \theta = \tilde{\theta} + \theta'.
\end{align}

Periodic boundary conditions are imposed on the lateral boundaries, whereas slip
boundary conditions are assumed elsewhere.  The Euler equations with the
prescribed initial and boundary data are integrated until $3,000$s of the model
time in order to obtain a solution for comparison.  If not mentioned otherwise
the simulations are carried out on a twice refined equally spaced grid resulting
in a spatial resolution of around $\Delta x=940$m.

\begin{figure}
  \begin{center}
  \includegraphics[width=0.95\textwidth]{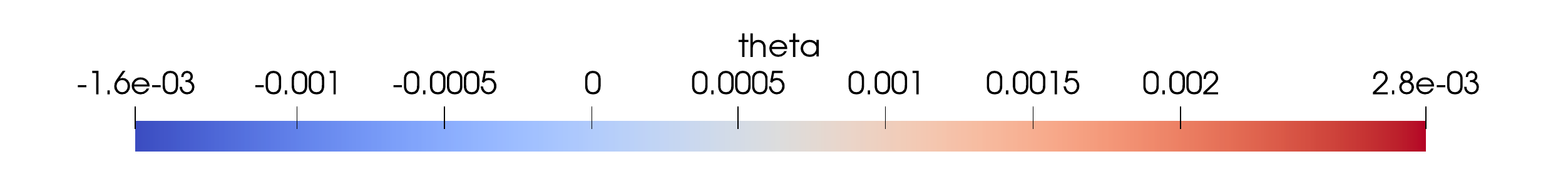} \\
  \subfloat[explicit $\Delta t = 0.028$]
  { \includegraphics[width=0.45\textwidth]{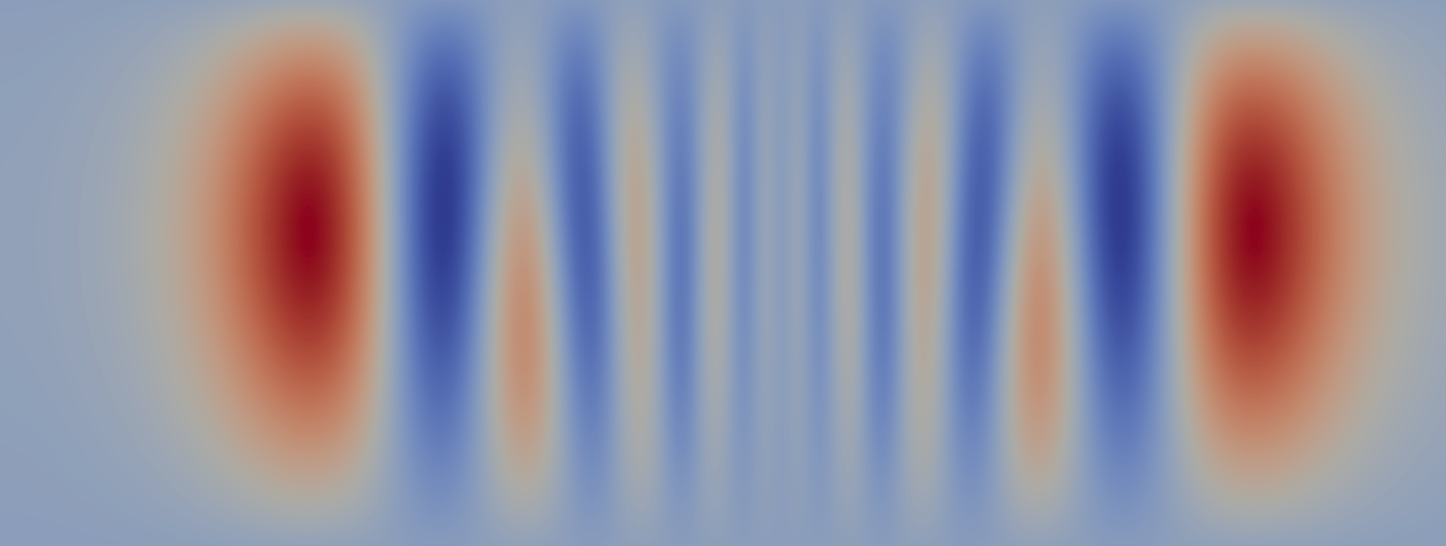}}
  \subfloat[implicit $\Delta t = 25$]
  { \includegraphics[width=0.45\textwidth]{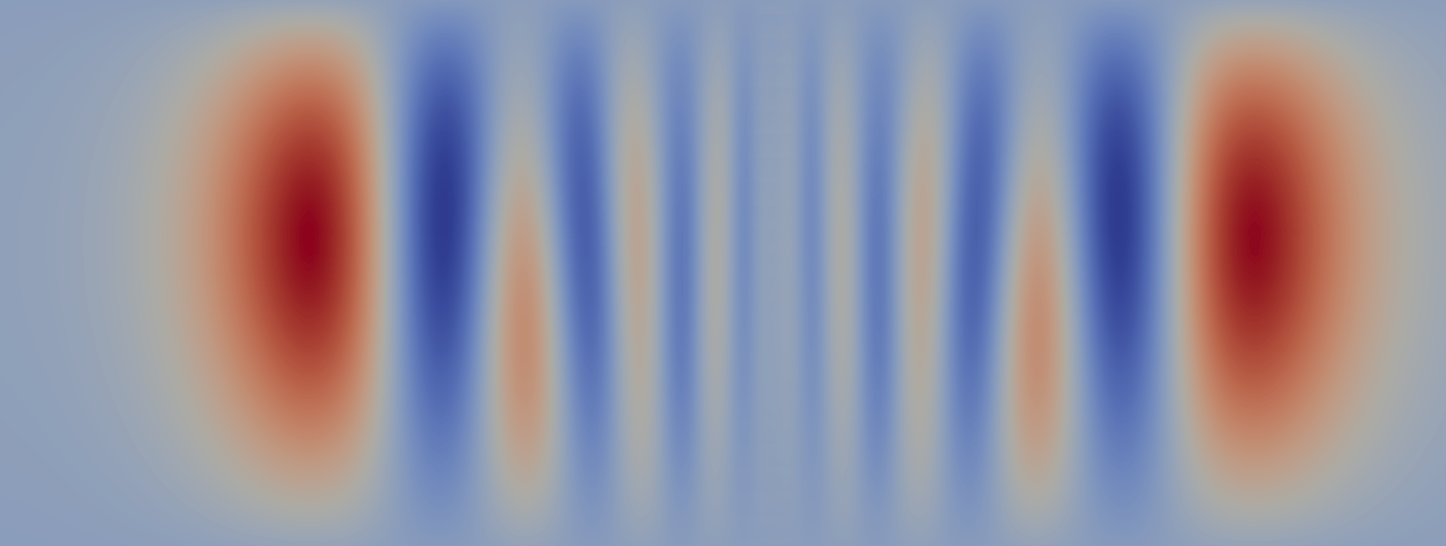}} \\
  \subfloat[implicit $\Delta t = 50$]
  { \includegraphics[width=0.45\textwidth]{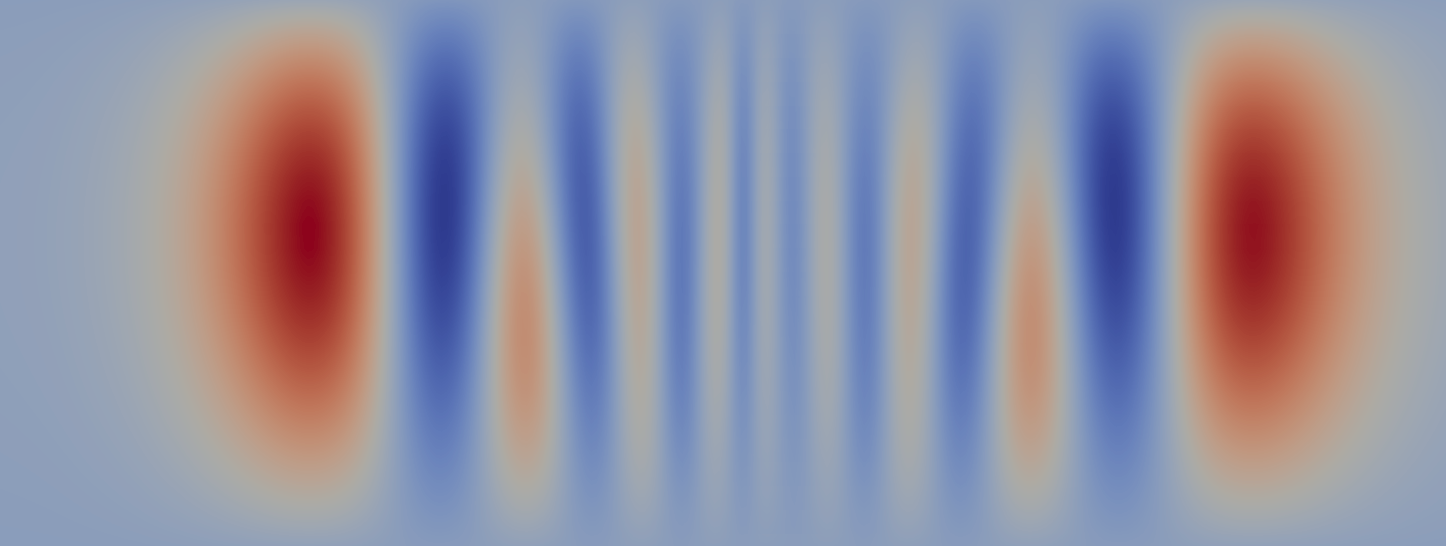}}
  \subfloat[implicit $\Delta t = 100$]
  { \includegraphics[width=0.45\textwidth]{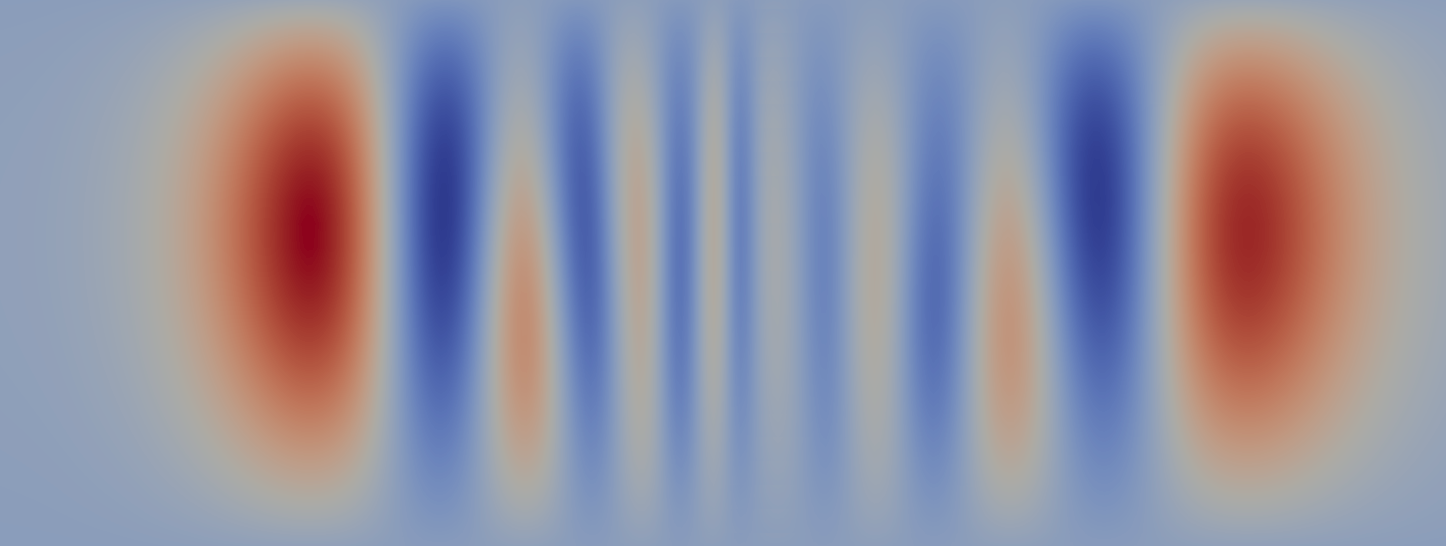}} \\
  \end{center}
  \caption{Inertia Gravity test case: perturbation in the potential temperature $\theta$ at the final time $T=3000$ for $\Delta x = 940$m
  using implicit solvers with different time steps and an explicit solver (top left).}
  \label{fig:igTheta}
\end{figure}

\begin{figure}
  \begin{center}
  \subfloat[comparison]
  { \includegraphics[width=0.45\textwidth]{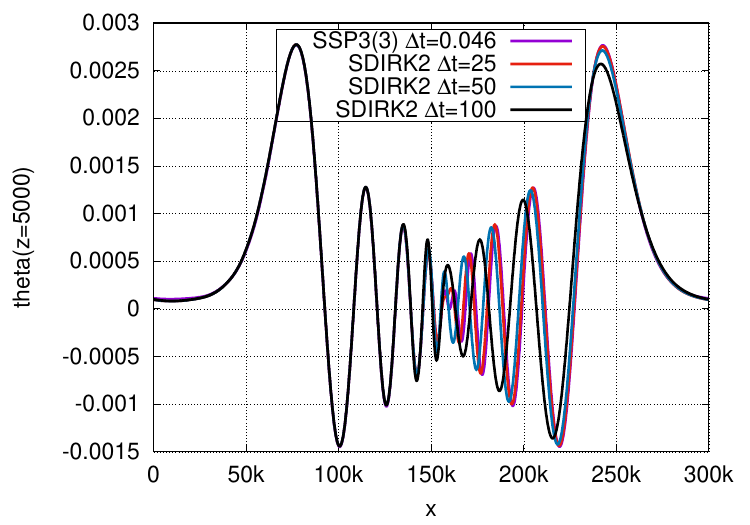}}
  \subfloat[zoomed comparison]
  { \includegraphics[width=0.45\textwidth]{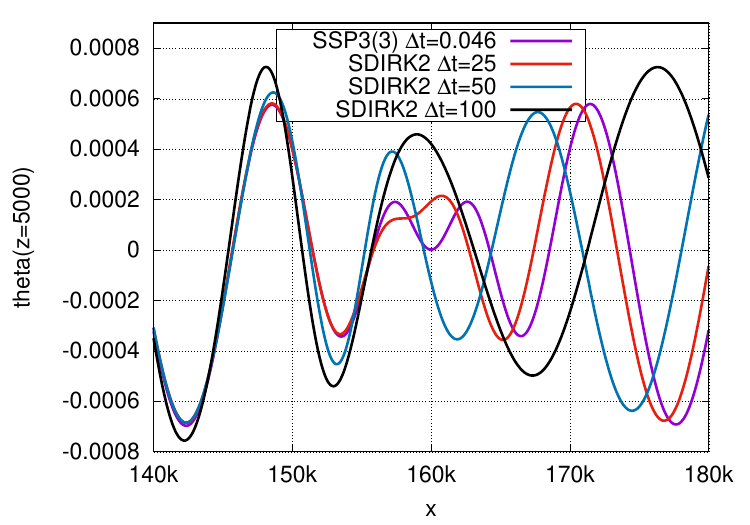}}
  \end{center}
  \caption{Inertia Gravity test case: perturbation in the potential temperature
  $\theta$ at the final time $T=3000$, $\Delta x = 940$m at $z=5,000$m.  The same
  time steps are used as in Figure~\ref{fig:igTheta}.  The solution with the
  explicit method in blue acts as a reference solution. The right figure shows
  $\theta$ around the center of the domain, results along the full range of $x$
  is shown on the left.  }
  \label{fig:igThetaCut1}
\end{figure}

\begin{figure}
  \begin{center}
  \subfloat[comparison]
  { \includegraphics[width=0.45\textwidth]{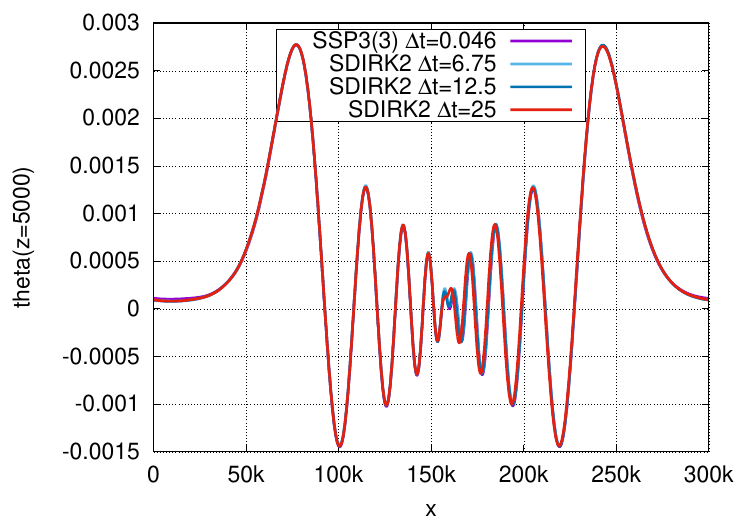}}
  \subfloat[zoomed comparison]
  { \includegraphics[width=0.45\textwidth]{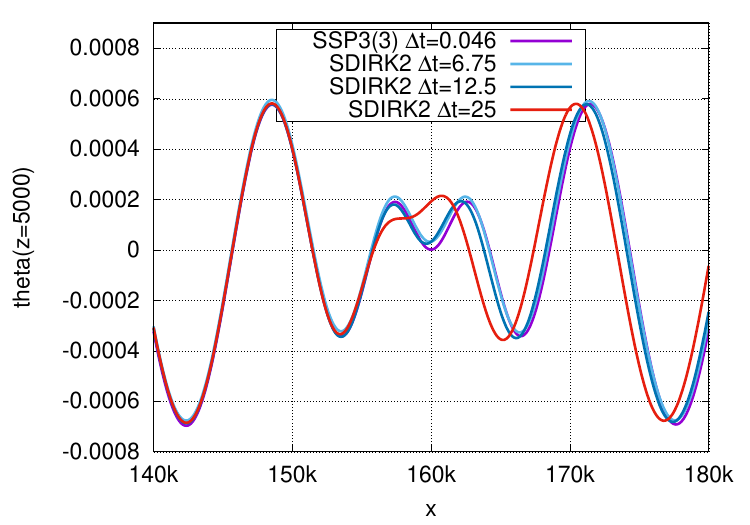}}
  \end{center}
  \caption{Inertia Gravity test case: perturbation in the potential temperature
  $\theta$ at the final time $T=3000$, $\Delta x = 940$m at $z=5,000$m.
  Smaller time steps are shown compared to Figure~\ref{fig:igThetaCut1}.}
  \label{fig:igThetaCut2}
\end{figure}

Figure~\ref{fig:igTheta} shows the potential temperature perturbation at the
final time using the implicit method with time steps $\Delta t=25,50,100$
together with a reference solution obtained with the explicit method using
$\Delta t=0.028$. The implicit methods are thus using time steps which are
around $1,000$ to $4,000$ times the stable time step of the explicit method. The
results looks still fairly similar although in the center some detail is lost
when using very large time steps.  This is clearly visible in
Figure~\ref{fig:igThetaCut1} where the solutions are plotted along the height
$y=5,000$. While the results with $\Delta t=25$ still broadly follow the
reference solution there is some clearer deviation when even larger time are
being used. Figure~\ref{fig:igThetaCut2} shows the same cut using time steps
$\Delta t=6.75,12.5,25$. It is clear that the results with $\Delta t=12.5$ which
is about $500$ times the time step used in the explicit simulation are very
close to the reference solution.

\begin{figure}
  \begin{center}
  \subfloat[$\Delta x=940m$, different $\Delta t$]
  {\includegraphics[width=0.49\textwidth]{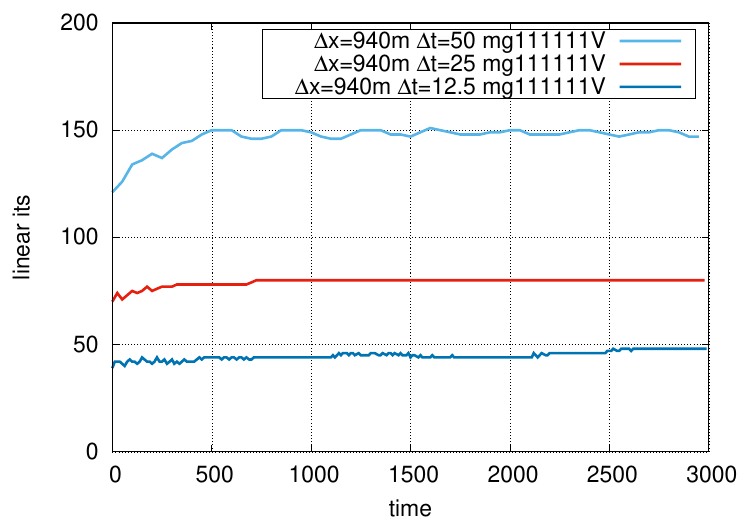}}
  \subfloat[at $\Delta x=470m$, different $\Delta t$, DG smoother on/off]
  { \includegraphics[width=0.49\textwidth]{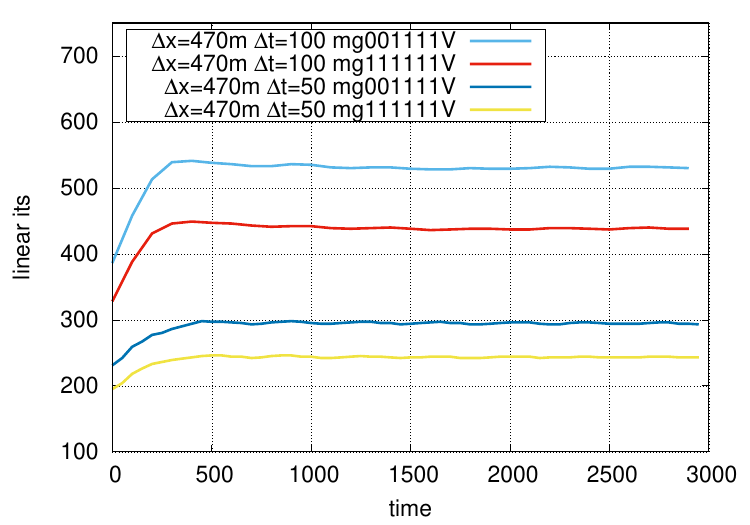}} \\
  \subfloat[at $\Delta x=470m$: DG smoother on/off, no precon]
  {
    \label{fig:igLinIts2}
    \includegraphics[width=0.49\textwidth]{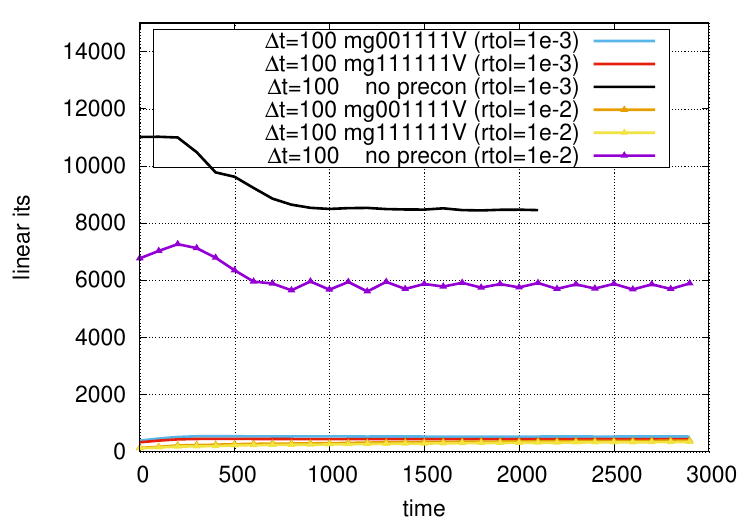}}
  \subfloat[$\Delta x=470m$, zoom of (c)]
  {\includegraphics[width=0.49\textwidth]{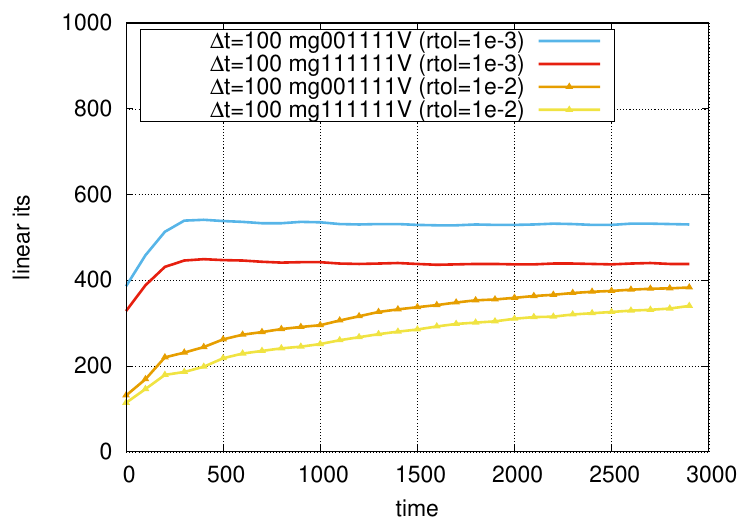}}

  \caption{Inertia Gravity test case: plots show the cumulative number of linear iterations needed to solve the non linear system
  in each time step over simulation time. Results with different values for the time step $\Delta t$ are shown
  with the plot on the left showing results for $\Delta x = 940m$ while the other two plots show results with $\Delta x = 470m$.
  Different choices of the multigrid parameters are compared including on the right results without preconditioning.}
  \label{fig:igLinIts}
  \end{center}
\end{figure}

To investigate the efficiency of the multigrid preconditioner, we show the
number of iterations needed in each time step in Figure \ref{fig:igLinIts}. The
figure on the left are for the same resolution ($\Delta x=940$m) used in the
previous plots. One can see that the number of iterations scales sublinearly
with increase in $\Delta t$ between $\Delta t=12.5$ and $\Delta t=25$ while
approximately twice the number of iterations are needed when going from $\Delta
t=25$ to $\Delta t=50$ indicating that $\Delta t$ around $25$ is a good choice
for this test case, i.e., about $500$ times the explicit time step. In the
center plot of Figure \ref{fig:igLinIts} we reduced the spatial resolution to
around $470$m by increasing the level to $\Delta x=470$m. We used a time step of
$\Delta t=50$ and in addition $\Delta t=100$ with two different parameter sets
for the multigrid method. Again a below linear increase in the number of
iterations can be seen. In addition the results indicate that the number of
iterations can be reduced by up to $20$\% by including pre and post smoothing
steps, although it should be noted that these do increase the computational cost
of each iteration step. The more detailed discussion of the choice of the
multigrid parameters is carried out for the next test case.  Finally, in Figure
\ref{fig:igLinIts2} we compare the implicit method with and without our
multigrid preconditioner. The number of iterations is significantly reduced by
our method to the point that with a smaller tolerance in the Newton solver the
simulation without preconditioning exceeded the time slot available on the
supercomputer.


%% file: sections/test-cases/rising-bubble.tex
The setup for the evolution of a warm bubble in a neutrally stratified atmosphere follows the description in \cite{robert:93}. The bubble, positioned at a certain height,
will start to rise and develop vortex structures in interaction with the cooler surrounding air. The stiffness in this test case comes from the difference in speed between the rising bubble and the fast sound waves
also present in the solution.

The governing equations are the compressible Euler equations with gravitational force. The domain is $[0, \ 1000] \times [0, \ 2000] $m$^2$ and the environmental atmosphere is determined by the formulas
\begin{eqnarray}
    \tilde{\theta} = 303.15 K, \quad \tilde{T} = T_0 - z g c_p^{-1}, \quad \tilde{p} = p_0 (\tilde{T}T_0^{-1})^{c_p/R_d}.
\end{eqnarray}
In the test case, $c_p = 1,005$J/(kg K), $c_v = 717.95$J/(kg K) and $g=9.80665$m/s$^2$. The initial atmosphere
is the sum of this neutrally stratified environmental atmosphere and a perturbation
induced by a deviation of the potential temperature
with the shape of a "flattened" Gaussian pulse, given by
\begin{eqnarray}
    \theta'(x, z) = A_0 \begin{cases}
        1.0, & r := ||(x, z) - (x_0, z_0)|| < a, \\
        exp(-(r-a)^2/s^2), &  0 \le r - a \le 3 s,\\
        0, & \text{else},
    \end{cases}
\end{eqnarray}
where $A_0 = 0.5$K, $x_0 = 500$m, $z_0 = 520$m, $a = 50$m and $s=100$m.
The introduction of the perturbation $\theta'$ into the initial state follows the approach described in \cite{dunecosmo:12}.
On all sides of the domain we impose slip boundary conditions. The system is integrated until 1200s of model time.

\begin{figure}[!ht]
  \subfloat[$t=0$]{\includegraphics[width=0.45\textwidth]{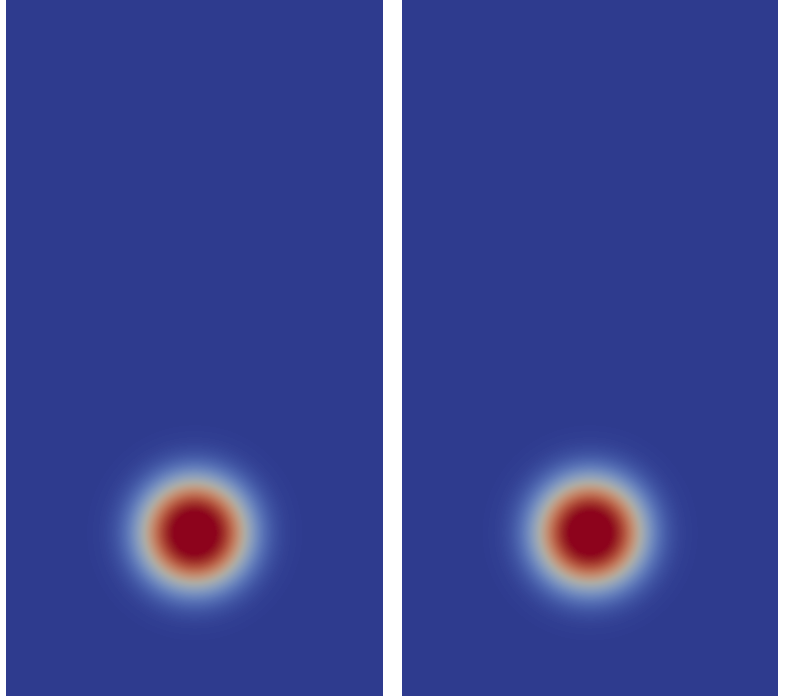}}
  \hfill
  \subfloat[$t=360$]{\includegraphics[width=0.45\textwidth]{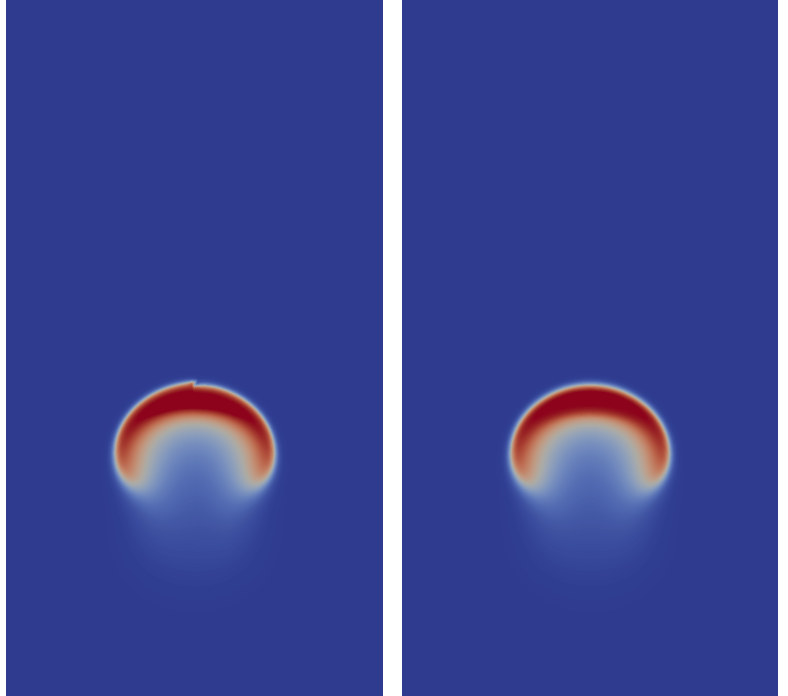}}\\
  \subfloat[$t=720$]{\includegraphics[width=0.45\textwidth]{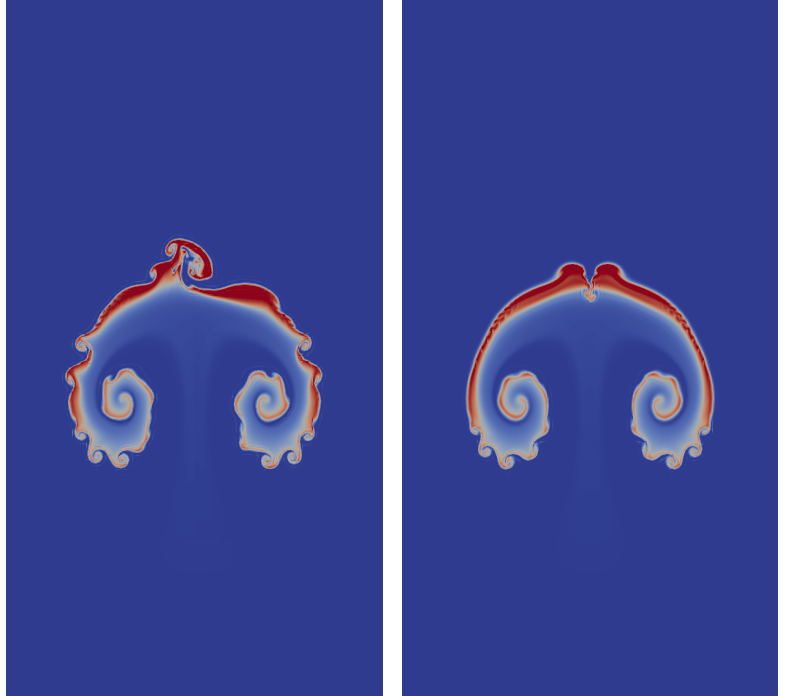}}
  \hfill
  \subfloat[$t=1200$]{\includegraphics[width=0.45\textwidth]{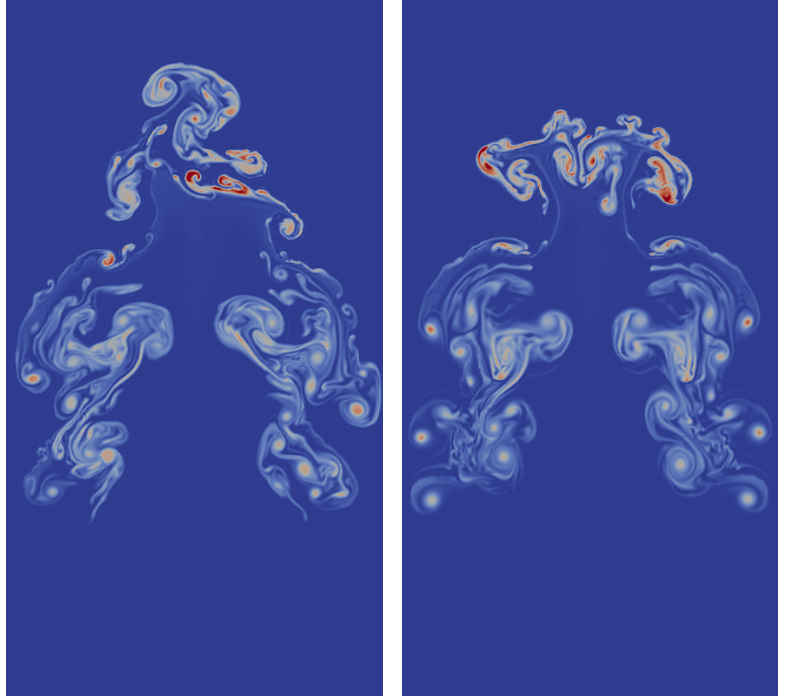}} \\
  \caption{Rising bubble test case: potential temperature perturbation at various times. Left, the solution
  computed at $\Delta x = 12.5m$ using the explicit time stepping with $\Delta t=0.002558s$, and right
  the solution at $\Delta x = 12.5m$ using implicit time stepping with multigrid preconditioning and time step $\Delta t = 5s$.
  Note that the implicit time stepping allows to use a time step that is $\approx 2000$ times larger than the explicit time step.
  Both schemes produce an acceptable solution.
  }
  \label{fig:risingbubble}
\end{figure}

For the runs at $\Delta x=25m$ we the DG solution is computed on level $l=2$ and for $\Delta x=12.5m$ we have $l=3$, which means $4$ or $5$ multigrid levels, respectively. For the simulations using the implicit SDIRK2 method, we use a fixed time step
of $\Delta t=5s$ and $\Delta t=10s$. Other time step sizes have been
tested but did yielded inferior results. The stable time step for the explicit method with a vertical resolution of $\Delta x = 12.5m$ was $\Delta t = 0.002558s$.
The solution of the test case at various points in time for the explicit and one of the considered implicit configurations is displayed in Figure \ref{fig:risingbubble}.

\begin{figure}[!ht]
  \begin{center}
    \setlength{\captionmargin}{0.5cm}
    \subfloat[Comparison of linear iteration counts for  $\Delta t=5s$ and $\Delta t=10s$]{\label{fig:linits_a}
    \includegraphics[width=0.45\textwidth]{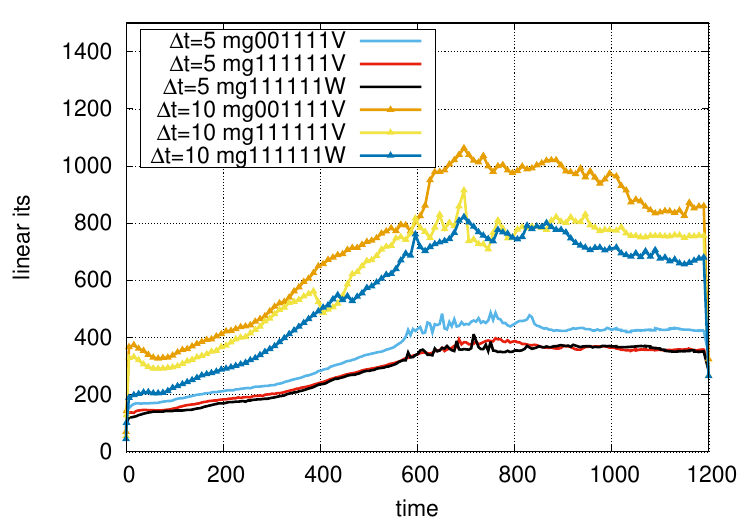}}
  \subfloat[Comparison of operator calls for $\Delta t=5s$ and $\Delta t=10s$]{\label{fig:linits_b}
    \includegraphics[width=0.45\textwidth]{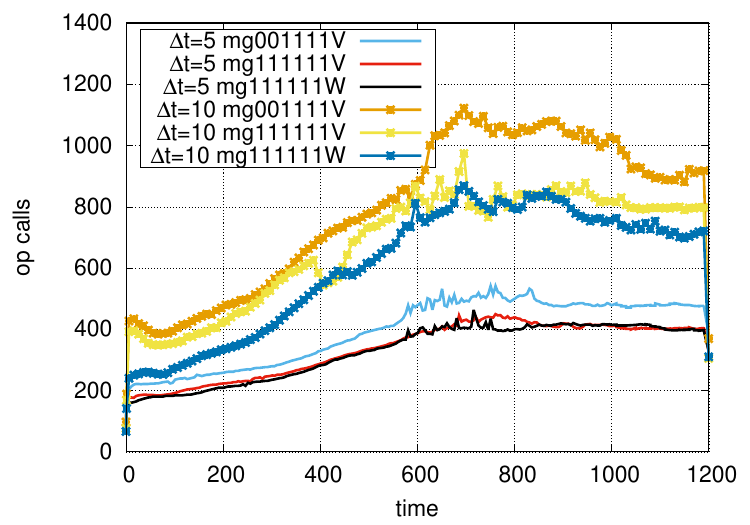}} \\
  \subfloat[Comparison of linear iteration counts for  $\Delta t=5s$ for two
    different grid width $\Delta x=25m$ and $\Delta x=12.5m$.]{\label{fig:linits_c}
    \includegraphics[width=0.45\textwidth]{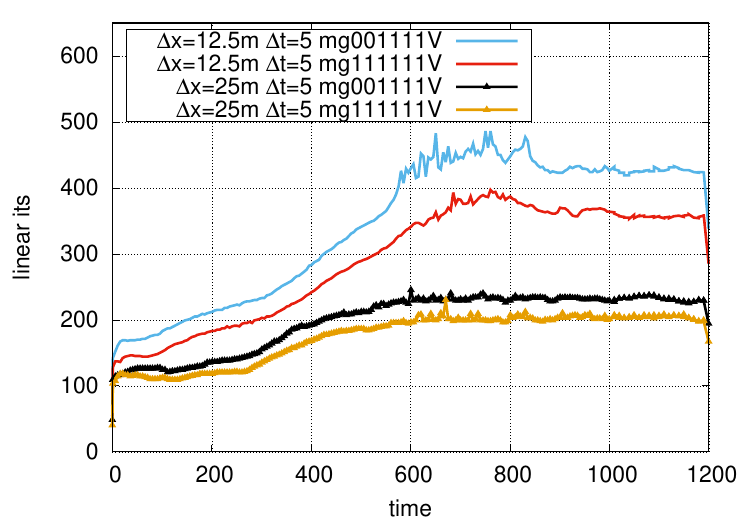}}
    \subfloat[Comparison including non-preconditioned run.]{\label{fig:linits_d}
    \includegraphics[width=0.45\textwidth]{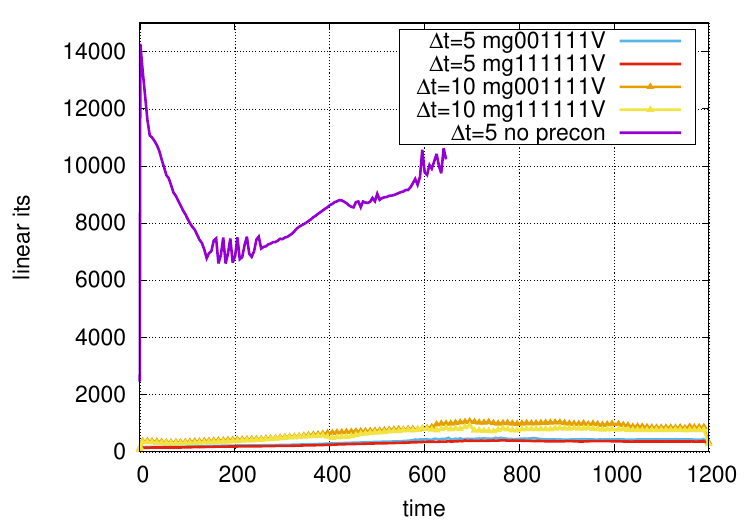}}
  \caption{Rising bubble test case: linear iterations for various configurations of the multigrid solver.}
  \label{fig:linits}
  \end{center}
\end{figure}

Figure \ref{fig:linits} contains a selection of different configurations of the
multigrid method for two time step sizes. In Figure \ref{fig:linits_a} we
compare the non-preconditioned implicit method with the preconditioned method.
The different configurations include either no pre
and post smoothing step on the DG solution and a mix of both, one or two smoothing steps on the FV levels as well as using a V-cycle and a W-cycle.
The two best configurations in terms of minimal linear iterations and minimal CPU time are $\text{\bf mg}001111V$ and $\text{\bf mg}111111V$,
where the additional pre and post smoothing on the DG solution reduces the number of linear iterations at the
cost of more operator evaluations.
In addition, a W-cycle configuration was tested, e.g. $\text{\bf mg}111111W$, yielding a competitive
number of linear iterations but the overall application of the multigrid method is much more expensive leading to doubling
in run time.


In Figure \ref{fig:linits_b} we compare the number of operator calls (DG and FV
combined) for the different combinations from Figure \ref{fig:linits_a}.
We can observe that the number of operator calls follow closely the number of linear
iterations. Therefore, the number of operator calls is not further considered in the following.

Figure \ref{fig:linits_c} compares the two best configurations from Figure
\ref{fig:linits_a} and \ref{fig:linits_b} for different  grid width $\Delta
x = 25m$ and $\Delta x=12.5m$. The number of linear iterations is
not grid independent but grows with decreasing grid width. However, in this
particular case the growth is less than a factor of two which we consider acceptable given that the time step size was kept constant at $\Delta t = 5s$.

Finally, Figure \ref{fig:linits_d} shows the linear iterations needed when
using the non-preconditioned implicit version. The simulation terminated prematurely
 due to exceeding the allocated time window on the super computer.

In summary, at grid width $\Delta x=12.5m$ the fastest run was $\text{\bf mg}111111V$ with $7915s$ to finish followed by $\text{\bf mg}0011111V$ with $8166s$, both using $\Delta t=5s$.
Compared to the explicit SSP3(4) time stepping, which needed $14000s$, this is roughly twice as fast.
All other multigrid configurations were slower, including those using a larger
time step size.
At $\Delta x=25m$, using a time step size of $\Delta t=10s$ the configuration $\text{\bf mg}111111V$ needed $1181s$ followed
by $\text{\bf mg}0011111V$ which $1196s$. The explicit  SSP3(4) time stepping took $1763s$.

%% file: sections/test-cases/density-current.tex
In \cite{straka:93}, the density current test case was proposed.
In this test case, one examines the evolution of a cold bubble in neutrally stratified atmosphere.
The bubble is positioned at a certain height, and it will start to fall and eventually hit the ground.
It will then slide along the ground level and create Kelvin-Helmholtz vortices.
Approximately every $300$s, a new Kelvin-Helmholtz vortex should appear.
In order to obtain grid-convergent solutions, enough viscosity should be employed.
In \cite{dunecosmo:12}, the authors have discussed that setting $\mu = 75$m$^2$/s is enough
to obtain grid convergence. The setup of the initial conditions is as follows.
The domain is $[0,\ 25,600]\times[0,\ 6,400]$m$^2$. The environmental atmosphere is determined by the formulas
\begin{align}
   \tilde{\theta} = 300 K, \quad \tilde{T} = T_0 - z g c_p^{-1}, \quad \tilde{p} = p_0 (\tilde{T}T_0^{-1})^{c_p/R_d}.
\end{align}
For this test case, we choose $c_p = 1,004$J/(kg K), $c_v=717$J/(kg K), and $g = 9.81$m/s$^2$ as in \cite{straka:93, dunecosmo:12}.
The initial data is a sum of this neutrally stratified environmental atmosphere and a
perturbation field induced by the deviation of the potential temperature
\begin{align}
    \theta'(x, z) = \begin{cases}
        \frac{\theta_c}{2}\big(1 + \cos(\pi r (x, z))\big), & \text{for}\ r(x, z) < 1, \\
        0, & \text{else},
    \end{cases}
\end{align}
where $\theta_c = -15$K is the maximal deviation of the potential
temperature from the stratified atmosphere, and $r$ is the distance to the center of the perturbation bubble, given by
\begin{align}
    r^2(x, z) = ((x-x_c)/x_r)^2 + ((z-z_c)/z_r)^2.
\end{align}

\begin{figure}[!ht]
  \begin{center}
  \includegraphics[width=0.95\textwidth]{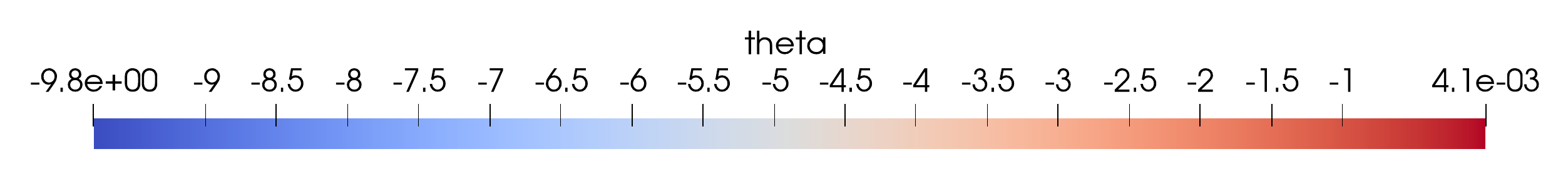} \\
  \subfloat[explicit $\Delta t = 0.01$]
  { \includegraphics[width=0.45\textwidth]{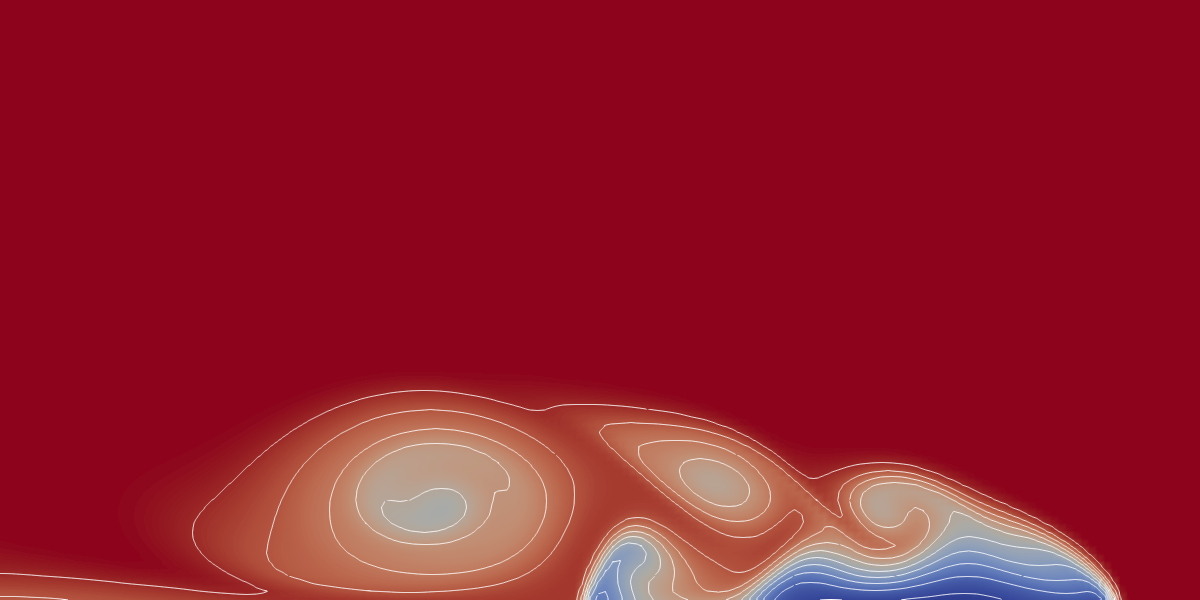}}
  \subfloat[implicit $\Delta t = 1$]
  { \includegraphics[width=0.45\textwidth]{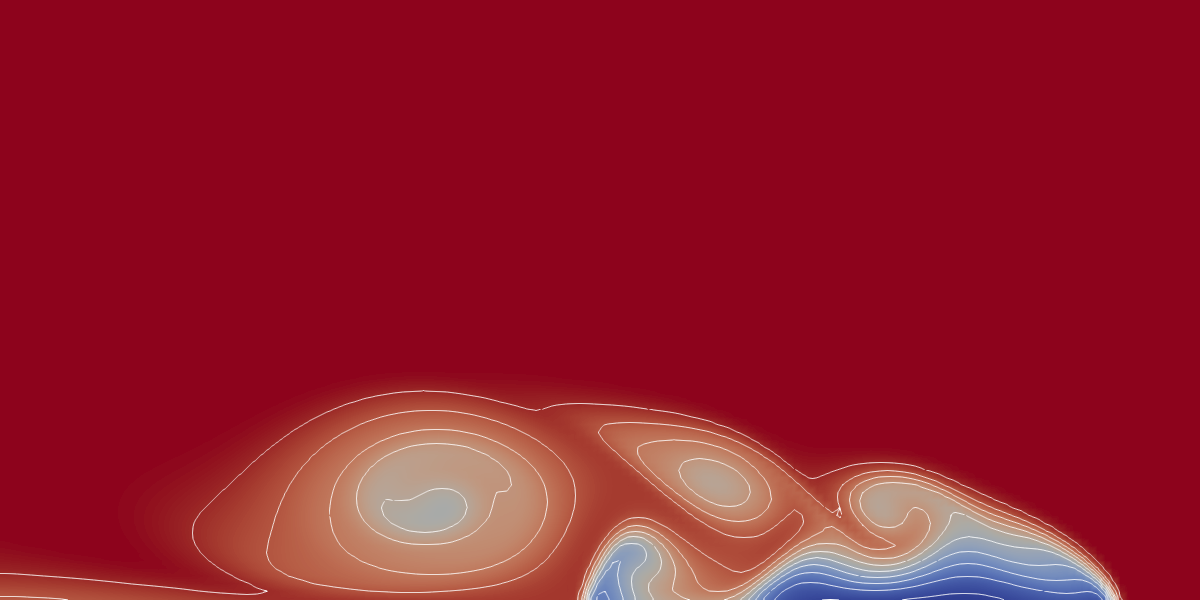}} \\
  \subfloat[implicit $\Delta t = 3$]
  { \includegraphics[width=0.45\textwidth]{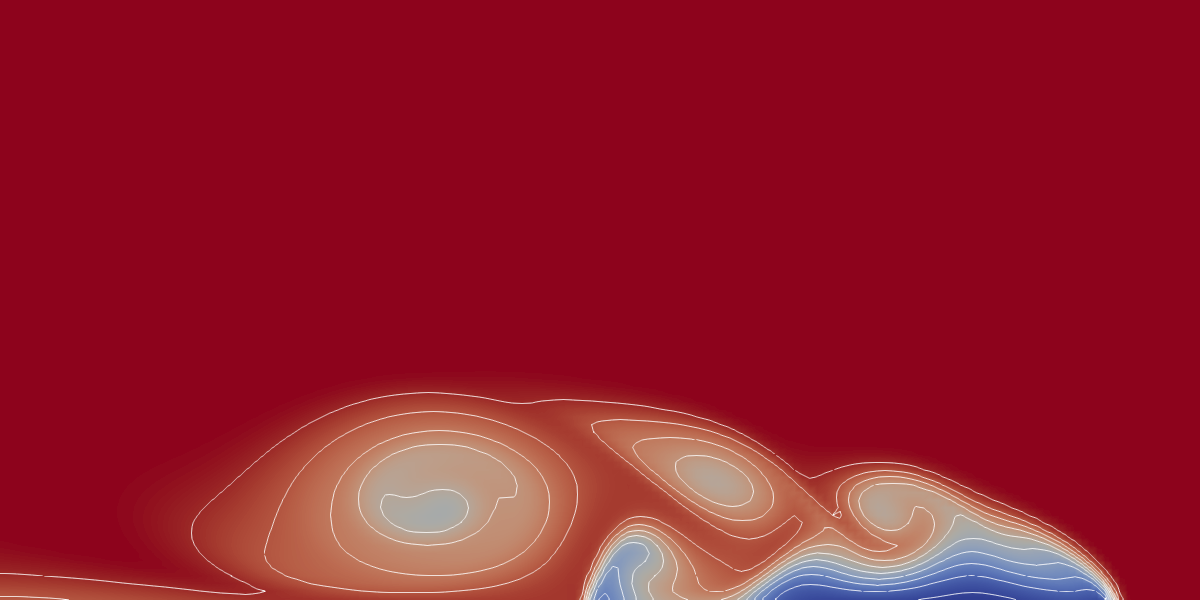}}
  \subfloat[implicit $\Delta t = 5$]
  { \includegraphics[width=0.45\textwidth]{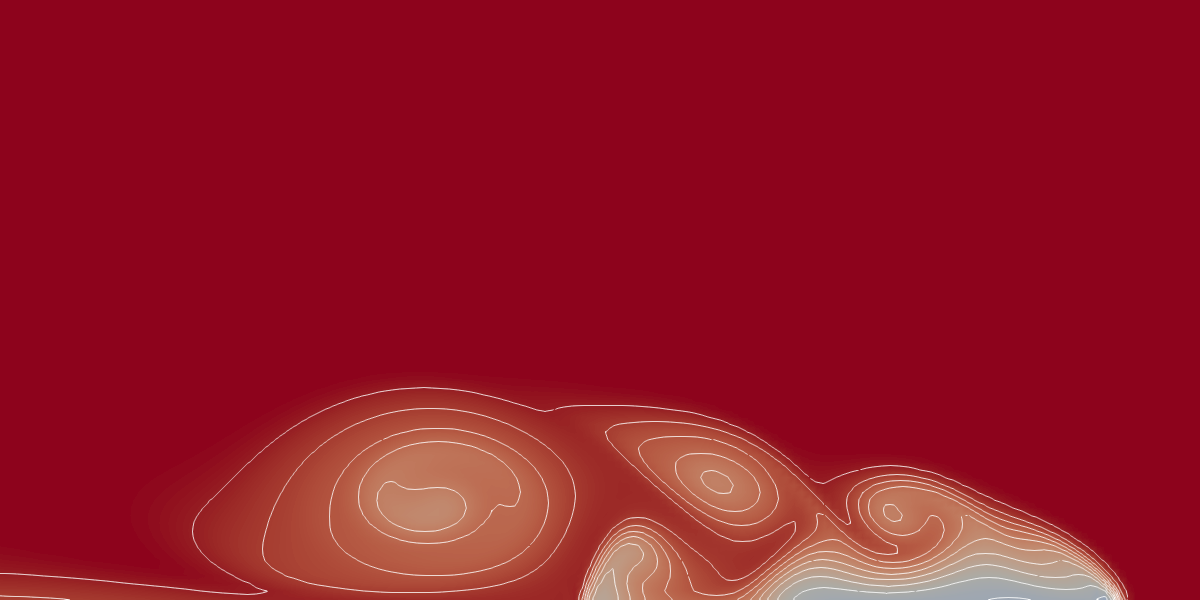}}
  \end{center}
  \caption{Density current test case: perturbation in the potential temperature
  $\theta$ at the final time $T=900$ for $\Delta x = 160$m using implicit
  solvers with different time steps and an explicit solver (top left). The
  solution patterns corresponds to those reported in \cite{dunecosmo:12}.  For
  $\Delta t = 5$ a deviation from the expected pattern is observed. Compare also
  with Figure \ref{fig:dcCut}.  }
  \label{fig:strakaTheta}
\end{figure}

\begin{figure}[!ht]
  \begin{center}
  \includegraphics[width=0.6\textwidth]{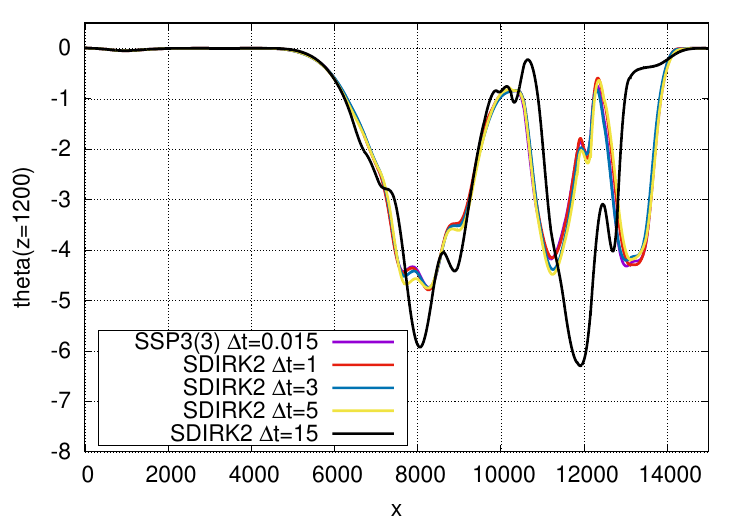}
  \end{center}
  \caption{Density current test case: solutions along the x-axis at $z=1200$m using
  different values of $\Delta t$. The solution using an explicit method is provided as reference solution.
  A good agreement can be observed up to and including $\Delta t=3$.
  }
  \label{fig:dcCut}
\end{figure}

\begin{figure}[!ht]
  \begin{center}
  \subfloat[$\Delta x=320$m, different $\Delta t$]
  { \includegraphics[width=0.49\textwidth]{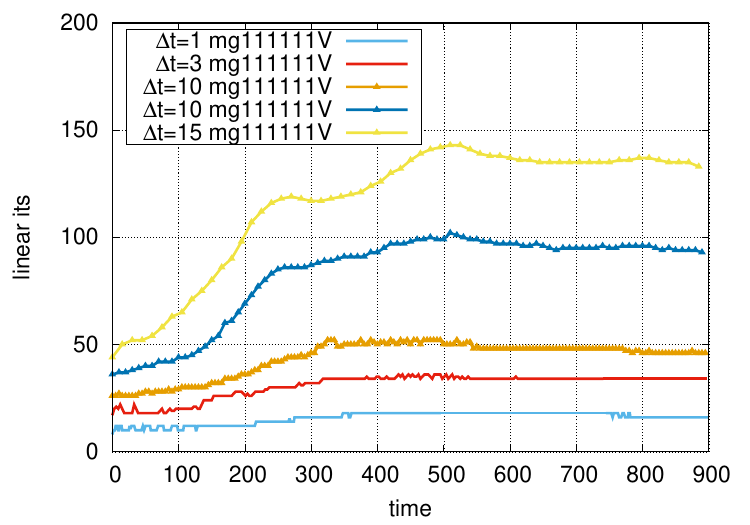} }
  \subfloat[$\Delta x=160$m, different $\Delta t$]
  { \includegraphics[width=0.49\textwidth]{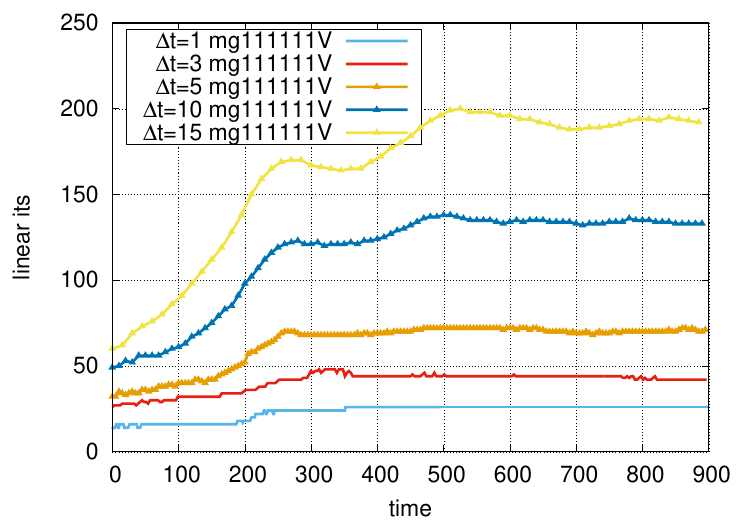} }
  \end{center}
  \caption{Density current test case: plots show the cumulative number of linear iterations needed to solve the non linear system
  in each time step over simulation time. Results with different values for the
  time step $\Delta t$ are shown using $\Delta x=320$m (left) and $\Delta
  x=160$m (right).}
  \label{fig:dcIters}
\end{figure}

The center of the perturbation bubble for the initial data is at $(x_c, z_c) = (0,\ 3,000)$m,
with radii $x_r = 4,000$m and $z_r = 2,000$m. The introduction of the perturbation $\theta'$
into the initial state follows the same way as in Eq. \eqref{eq:add-perturbation}.

On all sides of the domain we impose slip boundary conditions
(in spite of the fact that diffusion terms are present). The system is integrated until $900$s of the model time.

In Figure~\ref{fig:dcCut} we can compare the performance of the implicit solver
using different values for the time step $\Delta t$ compared with a reference
solution computed using the explicit method. Clearly using a $300$ times larger
time step compared to the stable time step needed in the explicit setting still
results in a accurate solution. Close to a factor of $1,000$ the solution shows
some clear deviation not accurately resolving all the vortices. We conclude with
again comparing the number of iterations needed with implicit solver using
different time steps and different spatial resolutions in
Figure~\ref{fig:dcIters}. Results are similar to what we already observed with a
sublinear growth in the number of iterations compared to increasing $\Delta t$.
These results also indicate that for this test case there is only an increase of
less then $50$\% in the number of iterations when the spatial resolution in
increased. This is even less of an increase compared to the rising bubble test
case.

At the finer resolution ($\Delta x=160$m) the explicit time stepping scheme needs about $2734$s to
  compute the solution at $T=900$.
The implicit scheme using the configuration $\text{\bf mg}111111V$ with $\Delta t = 3$ needed $1366$s and
$2580$s for $\Delta t=1$. At the coarser resolution ($\Delta x=320$m) the implicit scheme was slower
than the explicit scheme for all computed time step sizes.

%% file: sections/conclusions.tex
\section{Summary and Conclusions}

\label{sec:conclusion}

We considered implicit time integration for parallel high order DG methods for
compressible fluid flows. The arising systems are solved using a Jacobian-free
Newton-GMRES solver with preconditioner. The preconditioner is based on a
multigrid method for a first order finite volume discretization on an auxiliary
mesh, based on the DG mesh. To transfer between the meshes in a computationally
efficient manner, either a mass conservative $L_2$ mapping, or interpolation
with a mass fix can be used. As smoothers, low-storage explicit Runge-Kutta
pseudo time iterations are used, which are implemented in parallel in a
Jacobian-free manner. The preconditioner is thus Jacobian-free, uses little
extra memory, and achieves close to optimal arithmetic intensity for large
problems on parallel machines.

The numerical experiments are performed in the software framework DUNE-FEM for
three well known atmospheric test cases in 2D. These are a non-hydrostatic
inertia gravity case, a rising bubble of warm air, and a density current test
case. Gravity is incorporated in a well-balanced manner. The implicit high order
solver produces accurate results even for large CFL numbers. For these, the
unpreconditioned solver is prohibitively slow or does not converge at all. The
preconditioner increases the stability of the solver while substantially
reducing the number of iterations, making the implicit method competitive.
On grids with higher resolution the preconditioned implicit solver becomes
competitive in terms of run times compared to an explicit solver.

%% file: sections/appendix.tex
\section{Modified Newton-Cotes formula for finite volume cell centers}
\label{sec:quadrature}

A quadrature exact for polynomials of degree $k=3$ with quadrature points corresponding to
the cell centers of the finite volume grid from Figure \ref{fig:lgl_fv_points}
is given by the following quadrature points $G_i$ and weights $\omega_i$ for $i=0,1,2,3$ in one dimension 
\begin{align*}
  G_0 &= \frac{1}{8},\ \omega_0 = \frac{1625}{6000}  
  \quad G_1 = \frac{3}{8},\ \omega_1 = \frac{1375}{6000} \\
  G_2 &= \frac{5}{8},\ \omega_2 = \frac{1375}{6000} \quad  
  G_3 = \frac{7}{8},\ \omega_3 = \frac{1625}{6000} 
\end{align*}
Appropriate quadratures for higher dimensions are formed through tensor
product of the 1d quadrature. Weights for other polynomial degrees can be found
in the class \cpp{ModifiedNewtonCotes} in the \dunefem library
\cite{dune:Fem}.
